\documentclass[12pt,letterpaper,titlepage]{amsart}
\usepackage{amsmath, amssymb, amsthm, amsfonts,amscd,xr}
\newcommand{\N}{\mathbb{N}}
\newcommand{\Q}{\mathbb{Q}}
\newcommand{\Z}{\mathbb{Z}}
\newcommand{\R}{\mathbb{R}}
\newcommand{\C}{\mathbb{C}}

\def\beq{\begin{equation}}
\def\eeq{\end{equation}}

\def\a{\alpha}
\def\b{\beta}

\def\d{\delta}
\def\e{\varepsilon}

\def\A {\mathcal A}

\def\H{\mathcal H}
\def\O {\mathcal O}

\def\F{\mathbf F}

\newenvironment{res}
               {\begin{equation}
\begin{minipage}{0.85\textwidth}}
               { \end{minipage}\end{equation} }
\def\ber{\begin{res} }
\def\eer{\end{res}}

\numberwithin{equation}{section}
\newtheorem{thm}{Theorem}[section]

\newcommand{\ste}{\hfill\break}

\newtheorem{lemma}[thm]{Lemma}
\newtheorem{lem}[thm]{Lemma}
\newtheorem{lm}[thm]{Lemma}
\newtheorem{cor}[thm]{Corollary}

\newtheorem{prop}[thm]{Proposition}

\newtheorem{dfn}[thm]{Definition}
\newtheorem{rem}[thm]{Remark}

\makeatletter
\def\section{\@startsection {section}{1}{\z@}{3.5ex plus 1ex minus
    .2ex}{2.3ex plus .2ex}{\large\bf}}
    \def\subsection{\@startsection{subsection}{2}{\z@}{3.25ex plus 1ex minus
 .2ex}{1.5ex plus .2ex}{\bf}}
\makeatother

\def\qed{\hfill $\square$\par\vspace{5pt}}

\def\a{\alpha}
\def\b{\beta}
\def\l{\lambda}
\def\e{\epsilon}

\def\d{\delta}

\def\cc{C^\infty}
\def\cf{C^*_r(\H)}

\def\hf{\mathfrak{H}}

\def\nc{\mathcal{N}}

\def\q{{\bf q}}
\def\s{\sigma}

\def\A{\mathcal{A}}
\def\Ab{\mathbb{A}}

\def\O{\mathcal{O}}

\def\Ze{\mathcal{Z}}

\def\I{\mathcal{I}}
\def\A{\mathcal{A}}

\def\H{\mathcal{H}}
\def\F{\mathcal{F}}
\def\J{\mathcal{J}}
\def\K{\mathcal{K}}

\def\P{\mathcal{P}}
\def\S{\mathcal{S}}
\def\V{\mathcal{V}}
\def\W{\mathcal{W}}
\def\Ri{\mathcal{R}}

\def\waf{W^{\mathrm{aff}}}
\def\rnr{R_{\mathrm{nr}}}

\def\i{\iota}

\makeindex
\begin{document}
\title{The Schwartz algebra of an affine Hecke algebra}
\author{Patrick Delorme}
\address{Institut de Math\'ematiques de Luminy, UPR 9016 du
CNRS\\
Facult\'e des Sciences de Luminy\\
163 Avenue de Luminy, Case 901, 13288 Marseille Cedex 09\\
France\\
email: delorme@iml.univ-mrs.fr}
\author{Eric M. Opdam}
\address{Korteweg de Vries Institute for Mathematics\\
University of Amsterdam\\
Plantage Muidergracht 24\\
1018TV Amsterdam\\
The Netherlands\\
email: opdam@science.uva.nl}
\date{\today}
\thanks{During the preparation of this paper the second
named author was partially supported by a Pionier grant of
the Netherlands Organization for Scientific Research (NWO)}
\maketitle \tableofcontents
\section{Introduction}
An affine Hecke algebra is associated to a root datum (with basis)
$\Ri=(X, Y, R_{0}, {R}_{0}^\vee, F_{0})$, where $X, Y$ are
lattices with a perfect pairing, $R_{0}\subset X$ is a reduced
root system, ${R}^\vee_{0}\subset Y$ is the coroot system
and $F_{0}$ is
a basis of $R_{0}$, together with a length multiplicative
function $q$ of
the affine Weyl group associated to ${\mathcal R}$. It is denoted by
${\mathcal H}({\mathcal R}, q)$ or simply ${\H}$. It admits a natural
prehibertian structure (provided $q$ has values in $\R_+$, which
we assume throughout),
and it acts on its completion $L_{2}(\H)$
through bounded operators. Thus $\H$ is a Hilbert algebra in the
sense of \cite{Dix}.

The spectral decomposition of the left and right
representation of this Hibert algebra has been made explicit by
one of the authors (E.O., \cite{O}).
We will denote by $\F$ the isomorphism
between $L_2(\H)$ and its decomposition into irreducible representations.
We will call this map $\F$ the Fourier transform.

Another interesting completion of $\H$ is the
Schwartz algebra $\S\subset L_{2}(\H)$, which is a Fr\'echet
algebra completion of $\H$ \cite{O}.
The main theorem of this article is the
characterisation of the image of $\S$ by the Fourier transform
$\F$. This characterisation has several important consequences
which are described in Section \ref{sub:mainapp}. Let us
briefly discuss these applications.

First of all, we obtain the analog of Harish-Chandra's
completeness Theorem for generalized principal series
of real reductive groups. The representations involved
in the spectral decomposition of $L_{2}(\H)$ are, as
representations of $\H$, subrepresentation of certain finite
dimensional induced representations from parabolic subalgebras
(which are subalgebras of $\H$ which themselves belong to
the class of affine Hecke algebras).
We call these the standard tempered induced representations.
There exist standard interwining operators (see \cite{O})
between the standard induced tempered representations.
The completeness Theorem states that the commutant of the
standard tempered induced representations is generated by the
self-intertwining operators given by standard intertwining
operators.

Next we determine the image of the center of $\S$ and, as a
consequence,
we obtain the analog of Langlands' disjointness Theorem for real
reductive groups: two standard tempered induced representations
are either disjoint, i.e. without simple subquotient in common, or
equivalent.


Then we discuss the characterisation
of the Fourier transform, and of the set of minimal
central idempotents of the reduced $C^*$-algebra
$\mathcal{C}_r^*(\H)$ of $\H$.

Finally we observe that
that the dense subalgebra $\S\subset\mathcal{C}^*_r$
is closed for holomorphic calculus.

Let us now comment on the proof of the Main Theorem. As it
is familiar since Harish-Chandra's work on real reductive groups
\cite{HC}, \cite{HC2}, the determination of the image of
$\S$ by $\F$ requires a theory of the constant term
(see also \cite{DT} for the case of the hypergeometric Fourier
transform) for coefficients of tempered representations of
$\H$.  This theory is fairly simple using the decomposition of
these linear forms on $\H$ along weights of the action of the
abelian subalgebra $\A$ of $\H$. This subalgebra admits as a basis,
the family $\theta_{x}$, $x \in X$, which arises in the
Bernstein presentation of $\H$.

There is a natural
candidate ${\hat \S}$ for the image of $\S$ by $\F$. The inclusion
$\F(\S)\subset  {\hat \S} $ is easy to prove, using estimates of
the coefficients of standard induced tempered representations.

The only thing that remains to be proved at this point,
is that the inverse
of the Fourier transform, the wave packet operator $\J$,
maps ${\hat \S}$ to $\S$.  For this a particular role is played
by normalized smooth family of coefficients of standard tempered
induced  representations: these are smooth families divided by the
$c$-function. Of particular importance is the fact that the
constant terms of these families are finite sums of normalized
smooth families of coefficients for Hecke subalgebras of smaller
semisimple rank. This is a
nontrivial fact which requires the explicit computation of the
constant term of coefficients for generic standard tempered
induced representations. If $\I$ is the maximal ideal of the center
${\Ze}$ of $\H$ which annihilates such a representation, its
coefficients can be viewed as linear forms on Lusztig's
formal completion of $\H$ associated to $\I$. This allows to use
Lusztig's first reduction Theorem \cite{Lu}
which  decomposes this algebra. Some results on Weyl groups
are then needed to achieve this computation of the constant
term.

Once this property of normalized smooth family is
obtained, it easy to form wave packets in the Schwartz space, by
analogy with Harish-Chandra's work for real reductive groups \cite{HC}.
Simple lemmas on spectral projections of matrices and an
induction argument, allowed by the theory of the constant term,
lead to the desired result.

The paper is roughly structured as follows. First we discuss
in Sections 2 and 3 the necessary preliminary material on
the affine Hecke algebra and the Fourier transform on $L_2(\H)$.
We formulate the Main Theorem in Section 4, and we discuss
some of its consequences. In Section 5 we compute the constant
terms of coefficients of the standard induced representations and
of normalized smooth families of such coefficients. In Section 6
we use this and the material in the Appendix on spectral projections
in order to prove the Main Theorem. Finally, in the Appendix on the
$c$-function we have collected some fundamental properties of
the Macdonald $c$-functions on which many of our results
ultimately rely.
\section{The affine Hecke algebra and the Schwartz algebra}
This section serves as a reminder for the definition of the affine
Hecke algebra and related analytic structures. We refer the reader
to \cite{O}, \cite{Lu} and \cite{EO} for further background
material.
\subsection{The root datum and the affine Weyl group}
A reduced root datum is a $5$-tuple
$\Ri=(X,Y,R_0,R_0^\vee,F_0)$\index{R@$\Ri$, root datum},
where $X,Y$\index{X@$X,Y$, lattices}
are lattices with perfect pairing
$\langle\cdot,\cdot\rangle$\index{<@$\langle\cdot,\cdot\rangle$!a@perfect
pairing between $X$ and $Y$},
$R_0\subset X$\index{R2@$R_0\subset X$, reduced integral root system}
is a reduced root system,
$R_0^\vee\subset Y$\index{R2@$R_0^\vee\subset Y$, coroot system}
is the coroot system (which is
in bijection with $R_0$ via the map $\a\rightarrow\a^\vee$), and
$F_0\subset R_0$\index{F@$F_0\subset R_0$,
simple roots of $R_0$}
is a basis of simple roots of $R_0$.
The set $F_0$ determines a subset $R_{0,+}\subset R_0$ of
positive roots.

The Weyl group
$W_0=W(R_0)\subset\operatorname{GL}(X)$\index{W2@$W_0$, Weyl
group of $R_0$} of $R_0$ is
the group generated by the reflections $s_\a$
in the roots $\a\in R_0$.
The set
$S_0:=\{s_\a\mid\a\in F_0\}$\index{S@$S_0$,
simple reflections of $W_0$}
is called the set of simple reflections of $W_0$.
Then $(W_0,S_0)$ is a finite Coxeter group.

We define the affine Weyl group $W=W(\Ri)$ associated to $\Ri$
as the semidirect product
$W=W_0\ltimes X$\index{W@$W$, affine Weyl group}.
The lattice $X$ contains the root lattice
$Q$\index{Q@$Q$!$Q=Q(R_0)$, root lattice},
and the normal subgroup
$\waf:=W_0\ltimes Q\lhd W$\index{W2@$\waf=W_0\ltimes Q\subset W$}
is a Coxeter
group whose Dynkin diagram is given by the affine extension of
(each component of) the Dynkin diagram of $R_0^\vee$.
The affine root system of $\waf$ is given by
$R^{\mathrm{aff}}=
R_0^\vee\times\Z\subset Y\times\mathbb{Z}$\index{R5@$R^{\mathrm{aff}}$,
affine root system}.
Note that $W$ acts on $R^{\mathrm{aff}}$.

Let $R^{\mathrm{aff}}_+$\index{R5@$R^{\mathrm{aff}}_\pm$,
positive (negative) affine roots}
be the set of positive affine roots defined by
$R^{\mathrm{aff}}_+=\{(\a^\vee,n)\mid n>0,\mathrm{\ or\ }
n=0\mathrm{\ and\ }\a\in R_{0,+}\}$.
Let $F^{\mathrm{aff}}$ denote the corresponding set of
affine simple roots.
Observe that $F_0^\vee\subset F^{\mathrm{aff}}$\index{F@$F^{\mathrm{aff}}$,
affine simple roots}.
If $S^{\mathrm{aff}}$\index{S@$S^{\mathrm{aff}}$,
simple reflections of $\waf$}
denotes the associated set of
affine simple reflections, then
$(\waf,S^{\mathrm{aff}})$ is an affine Coxeter group.

In this paper we adhere to the convention
$\N=\{1,2,3,\cdots\}$\index{N@$\N=\{1,2,3,\cdots\}$} and
$\Z_+=\{0,1,2,\dots\}$\index{Z@$\Z_+=\{0,1,2,\dots\}$}.
We define the length function
$l:W\to\mathbb{Z}_+$\index{l1@$l$, length function on $W$}
on $W$ as usual, by means of the formula
$l(w):=|R^{\mathrm{aff}}_+\cap w^{-1}(R^{\mathrm{aff}}_-)|$.
Let
$\Omega\subset W$\index{0Y@$\Omega$, length $0$ elements in $W$}
denote the set $\{w\in W\mid l(w)=0\}$.
It is a subgroup of $W$, complementary to $\waf$. Therefore
$\Omega\simeq X/Q$ is a finitely generated Abelian subgroup
of $W$.

Let
$X^+\subset X$
denote the cone of dominant elements
$X^+=\{x\in X\mid
\forall\a\in R_{0,+}:\
\langle x,\a^\vee\rangle\geq 0\}$\index{X@$X^+$, cone of dominant
elements in $X$}.
Then $Z_X:=X^+\cap
X^-\subset X$\index{Z9@$Z_X$, length $0$ translations in $W$}
is a sublattice which is
central in $W$. In particular it follows that
$Z_X\subset\Omega$. The quotient $\Omega_f\simeq
\Omega/Z_X$\index{0Y@$\Omega_f=\Omega/Z_X$}
is a finite Abelian group which acts faithfully
on $S^{\mathrm{aff}}$ by means of diagram automorphisms.

We choose a basis $z_i$ of $Z_X$, and define
a norm
$\Vert\cdot\Vert$\index{$\Vert\cdot\Vert$, norm
on the vector space $\Q\otimes_\Z Z_X$}
on the rational vector space
$\Q\otimes_\Z Z_X$ by $\Vert\sum l_iz_i\Vert:=\sum |l_i|$.
We now define a norm
$\mathcal{N}$\index{N9@$\mathcal{N}$, norm function on $W$}
on $W$ by
\begin{equation}\label{eq:norm}
\mathcal{N}(w):=l(w)+\Vert w(0)^0\Vert,
\end{equation}
where $w(0)^0$ denotes the projection of $w(0)$ onto
$\Q\otimes_\Z Z_X$ along $\Q\otimes_\Z Q$.
The norm $\mathcal{N}$ plays an important role in this
paper. Observe that it satisfies
\begin{equation}
\mathcal{N}(ww^\prime)\leq \mathcal{N}(w)+\mathcal{N}(w^\prime),
\end{equation}
and that $\mathcal{N}(w)=0$ iff $w$ is an element of $\Omega$
of finite order.

We call $\Ri$ semisimple if $Z_X=0$.
\subsection{Standard parabolic subsystems}
A root subsystem $R^\prime\subset R_0$ is called parabolic
if $R^\prime=\mathbb{Q}R^\prime\cap R_0$. The Weyl group
$W_0$ acts on the collection of parabolic root subsystems.
Let
$\P$\index{P@$\P$, power set of $F_0$}
be the power set of $F_0$. With
$P\in\P$\index{P@$P\in\P$, subset of $F_0$}\index{Q@$Q$!$Q\in\P$,
subset of $F_0$}
we
associate a standard parabolic root subsystem
$R_P\subset R_0$\index{R3@$R_P\subset R_0$, parabolic subsystem,
root system of $\Ri_P$}
by $R_P:=\mathbb{Z}P\cap R_0$.
Every parabolic root subsystem is $W_0$-conjugate
to a standard parabolic subsystem.

We denote by
$W_P=W(R_P)\subset W_0$\index{W4@$W_P$, Weyl group of $R_P$,
parabolic subgroup $W_0$}
the Coxeter
subgroup of $W_0$ generated by the reflections
in $P$. We denote by
$W^P$\index{W4a@$W^P=W_0/W_P$, set of left cosets $wW_P$.
If $P\subset F_0$,
identified with shortest length representatives}
the set of shortest length
representatives of the left cosets $W_0/W_P$ of
$W_P\subset W_0$.

Given $P\in\P$ we define a sub root datum
$\Ri^P\subset \Ri$\index{R1@$\Ri^P$, root datum associated to $P$}
simply by $\Ri^P:=(X,Y,R_P,R_P^\vee,P)$. We also define a ``quotient
root datum'' $\Ri_P$ of $\Ri^P$ by
$\Ri_P=(X_P,Y_P,R_P,R_P^\vee,P)$\index{R1@$\Ri_P$,
semisimple root datum associated to $P$}
where
$X_P:=X/(X\cap(R_P^\vee)^\perp)$\index{X1a@$X_P\supset R_P$,
lattice of $\Ri_P$,
character lattice of $T_P$}
and
$Y_P=Y\cap \mathbb{Q}R_P^\vee$\index{Y@$Y_P\supset R_P^\vee$,
lattice of $\Ri_P$,
cocharacter lattice of $T_P$}.
The root datum $\Ri_P$
is semisimple.
\subsection{Label functions and root labels}
A positive real label function is a
length multiplicative function
$q:W\to\mathbb{R}_+$\index{q@$q$, $l$-multiplicative function on
$W$}.
This means that $q(ww^\prime)=q(w)q(w^\prime)$ whenever
$l(ww^\prime)=l(w)+l(w^\prime)$, and that $q(\omega)=1$
for all $\omega\in \Omega$.

Such a function
$q$ is uniquely determined by its restriction to
the set of affine simple reflections $S^{\mathrm{aff}}$.
By the
braid relations and the action of $\Omega_f$ on
$S^{\mathrm{aff}}$
it follows easily that $q(s)=q(s^\prime)$ whenever
$s,s^\prime\in S^{\mathrm{aff}}$ are $W$-conjugate.
Hence there exists
a unique $W$-invariant function
$a\to q_a$\index{q@$q_a$, affine root label}
on
$R^{\mathrm{aff}}$ such that $q_{a+1}=q(s_a)$ for all
simple affine roots
$a\in F^{\mathrm{aff}}$\index{a5@$a=(\alpha^\vee,k)$, affine root}.

We associate a possibly non-reduced root system
$\rnr$\index{R4@$\rnr$, non reduced root system}
with $\Ri$ by
\begin{equation}
\rnr:=R_0\cup\{2\a\mid \a^\vee\in R_0^\vee\cap 2Y\}.
\end{equation}
If $\a\in R_0$ then $2\a\in\rnr$ iff the affine roots
$a=\a^\vee$ and $a=\a^\vee+1$ are not $W$-conjugate.
Therefore we can also characterize the label
function $q$ on $W$ by means of the following extension
of the set of root labels
$q_{\a^\vee}$\index{q@$q_{\alpha^\vee}$, label for $\a^\vee\in\rnr^\vee$}
to arbitrary $\a\in\rnr$.
If $\a\in R_0$ with $2\a\in\rnr$, then we define
\begin{equation}
q_{\a^\vee/2}:=\frac{q_{\a^\vee+1}}{q_{\a^\vee}}.
\end{equation}
With these conventions we have for all $w\in W_0$
\begin{equation}
q(w)=\prod_{\alpha\in R_{\mathrm{nr},+}\cap
w^{-1}R_{\mathrm{nr,-}}} q_{\alpha^\vee}.
\end{equation}

We denote by $R_1\subset X$\index{R4@$R_1$, system of long roots in $\rnr$}
the reduced root system
\begin{equation}
R_{1}:=
\{\alpha\in R_\mathrm{nr}\mid
2\alpha\not\in R_\mathrm{nr}\}.
\end{equation}
\subsubsection{Restriction to parabolic subsystems}
Let $P\in\P$. Both the non-reduced root system associated
with $\Ri^P$ and the non-reduced root system associated
with $\Ri_P$ are equal to $R_{P,nr}:=\mathbb{Q}R_P\cap\rnr$.
We define a collection of root labels
$q_{P,\a^\vee}=q^P_{\a^\vee}$\index{q@$q_P$, label
function on $\Ri_P$ determined by the
restriction of $q$ to $R_{P,\mathrm{nr}}$}\index{q@$q^P$,
label function on $\Ri^P$ determined by the
restriction of $q$ to $R_{P,\mathrm{nr}}$}
for $\a\in R_{P,\mathrm{nr}}$ by restricting
the labels of $\rnr$ to $R_{P,\mathrm{nr}}\subset\rnr$.
Then $q_P$ denotes the length-multiplicative function
on $W(\Ri_P)$ associated with this label function on
$R_{P,\mathrm{nr}}$, and $q^P$ denotes the associated
length multiplicative function on $W(\Ri^P)$.
\subsection{The Iwahori-Hecke algebra}
Given a root datum $\Ri$ and a (positive real)
label function $q$ on the associated affine Weyl group $W$,
there exists a unique associative complex Hecke algebra
$\H=\H(\Ri,q)$\index{H@$\H$, affine Hecke algebra}
with $\mathbb{C}$-basis $N_w$\index{N9@$N_w$, basis elements of $\H$}
indexed by
$w\in W$, satisfying the relations
\begin{enumerate}
\item $N_{ww^\prime}=N_wN_{w^\prime}$
for all $w,w^\prime\in W$ such that
$l(ww^\prime)=l(w)+l(w^\prime)$.
\item
$(N_s+q(s)^{-1/2})(N_s-q(s)^{1/2})=0$
for all $s\in S^{\mathrm{aff}}$.
\end{enumerate}
Notice that the algebra $\H$ is unital, with unit $1=N_e$.
Notice also that it follows from the defining relations
that $N_w\in\H$ is invertible, for all $w\in W$.

By convention we assume that the label function $q$
is of the form
\begin{equation}\label{eq:conv}
\index{q@$\q$, base for the labels $q(s)$}
\index{f@$f_s=\log_{\q}(q(s))$}
q(s)=\q^{f_s}.
\end{equation}
The parameters $f_s\in\mathbb{R}$ are fixed, and
the base $\q$ satisfies $\q>1$.
\subsubsection{Bernstein presentation}
There is another, extremely important presentation of
the algebra $\H$, due to Joseph Bernstein (unpublished).
Since the length function is additive on the dominant cone
$X^+$, the map $X^+\ni x\to N_x$ is a homomorphism of the
commutative monoid $X^+$ with values in $\H^\times$, the group
of invertible elements of $\H$. Thus there exists a unique
extension to a homomorphism
$X\ni x\to \theta_x\in\H^\times$\index{0i@$\theta_x$,
basis elements of $\A$}
of the lattice $X$ with values in $\H^\times$.

The Abelian subalgebra of $\H$ generated by
$\theta_x$, $x\in X$, is denoted by
$\A$\index{A@$\A$, abelian subalgebra of $\H$}. Let
$\H_0=\H(W_0,q_0)$\index{H2@$\H_0=\H(W_0,q\mid_{S_0})\subset\H$}
be the finite type Hecke algebra
associated with $W_0$ and the restriction $q_0$ of $q$
to $W_0$. Then the Bernstein presentation asserts that
both the collections $\theta_xN_w$ and $N_w\theta_x$
($w\in W_0$, $x\in X$) are bases of $\H$, subject only
to the cross relation (for all $x\in X$ and
$s=s_\alpha$ with $\a\in F_0$):
\begin{gather}
\begin{split}
&\theta_xN_s-N_s\theta_{s(x)}=\\
&\left\{
\begin{array}{ccc}
&(q_{\alpha^\vee}^{1/2}-
q_{\alpha^\vee}^{-1/2})\frac{\theta_x-\theta_{s(x)}}
{1-\theta_{-\alpha}}&{\rm if}\ 2\alpha\not\in\rnr.\\
&((q_{\alpha^\vee/2}^{1/2}q_{\alpha^\vee}^{1/2}-
q_{\alpha^\vee/2}^{-1/2}q_{\alpha^\vee}^{-1/2})
+(q_{\alpha^\vee}^{1/2}-
q_{\alpha^\vee}^{-1/2})\theta_{-\alpha})
\frac{\theta_x-\theta_{s(x)}}
{1-\theta_{-2\alpha}}&{\rm if}\ 2\alpha\in\rnr.\\
\end{array}
\right.\\
\end{split}
\end{gather}
\subsubsection{The center $\Ze$ of $\H$}
A rather immediate consequence of the Bernstein
presentation of $\H$ is the description of the
center of $\H$:
\begin{thm}\label{thm:cent}
The center
$\Ze$\index{Z91@$\Ze$, the center of $\H$}
of $\H$ is equal to $\A^{W_0}$.
In particular, $\H$ is finitely generated
over its center.
\end{thm}
As an immediate consequence we see that irreducible
representations of $\H$ are finite dimensional by
application of (Dixmier's version of) Schur's lemma.

We denote by
$T$\index{T@$T=\operatorname{Hom}_\Z(X,\C^\times)$, complex algebraic torus}
the complex torus
$T=\operatorname{Hom}(X,\mathbb{C}^\times)$ of complex
characters of the
lattice $X$. The space $\operatorname{Spec}(\Ze)$
of complex homomorphisms of $\Ze$ is thus canonically
isomorphic to the (geometric) quotient
$W_0\backslash T$.

Thus, to an irreducible representation $(V,\pi)$
of $\H$ we attach an orbit $W_0t\in W_0\backslash T$,
called the central character
\index{central character of $(V,\pi)$}
of $\pi$.
\subsubsection{Parabolic subalgebras and their
semisimple quotients}
We consider another important consequence of the
Bernstein presentation of $\H$:
\begin{prop}\label{prop:parber}
\begin{enumerate}
\item
The Hecke algebra $\H^P:=\H(\Ri^P,q^P)$\index{H4@$\H^P=\H(\Ri^P,q^P)$,
parabolic subalgebra of $\H$}\index{H5@$\H_P=\H(\Ri_P,q_P)$,
semisimple quotient of $\H^P$}
is isomorphic to
the subalgebra of $\H$ generated by $\A$ and
the finite type Hecke subalgebra $\H(W_P):=\H(W_P,q|_{W_P})$.
\item
We can view $\H_P:=\H(\Ri_P,q_P)$ as a quotient of $\H^P$ via
the surjective homomorphism $\phi_1:\H^P\to\H_P$ characterized
by (1) $\phi_1$ is the identity on the finite type subalgebra
$\H(W_P)$ and (2) $\phi_1(\theta_x):=\theta_{\overline{x}}$,
where
$\overline{x}\in X_P$ is the canonical image of $x$ in
$X_P=X/(X\cap(R_P^\vee)^\perp)$.
\end{enumerate}
\end{prop}
Let $T^P$ denote the character torus of the lattice
$X/(X\cap\mathbb{Q}R_P)$. Then $T^P\subset T$ is a
subtorus which is fixed for all the elements $w\in W_P$
and which is inside the simultaneous
kernel of the $\a\in R_P$. The next result again
follows simply from the Bernstein presentation:
\begin{prop}\label{prop:partwist}
There exists a family of automorphisms $\psi_t$ ($t\in T^P$)
of $\H^P$,
defined by $\psi_{t}(\theta_x N_w)=x(t)\theta_{x} N_w$.
\end{prop}
We use the above family of automorphisms to twist the
projection $\phi_1:\H^P\to\H_P$.
Given $t\in T^P$, we define the epimorphism
$\phi_t:\H^P\to\H_P$ by
$\phi_t:=\phi_1\circ\psi_t$\index{0v@$\phi_{t^P}:\H^P\to\H_P$, family of
surjective homomorphisms}.
\subsection{Intertwining elements}
When $s=s_\alpha\in S_0$ (with $\alpha\in F_1$), we define an
intertwining element $\i_s\in\H$ as follows:
\begin{equation*}
\begin{split}
\i_s&=q_{\a^\vee}q_{2\a^\vee}(1-\theta_{-\alpha})N_s+
((1-q_{\alpha^\vee}q_{2\alpha^\vee})
+q_{\alpha^\vee}^{1/2}(1-q_{2\alpha^\vee})\theta_{-\alpha/2})\\
&=q_{\a^\vee}q_{2\a^\vee}
N_s(1-\theta_{\alpha})+((q_{\alpha^\vee}q_{2\alpha^\vee}-1)
\theta_{\alpha}
+q_{\alpha^\vee}^{1/2}(q_{2\alpha^\vee}-1)\theta_{\alpha/2})\\
\end{split}
\end{equation*}
\index{0j@$\iota_s$, intertwining element of $\H$}
(If $\a/2\not\in X$ then we put $q_{2\a^\vee}=1$;
see Remark \ref{rem:conv}.)
We recall from \cite{EO}, Theorem 2.8 that these elements satisfy the
braid relations, and they satisfy (for all $x\in X$)
\begin{equation}
\i_s\theta_x=\theta_{s(x)}\i_s
\end{equation}

Let
$\mathcal{Q}$\index{Q@$\mathcal{Q}$, quotient field of $\Ze$}
denote the quotient field of the centre
$\Ze$ of $\H$, and let
${}_\mathcal{Q}\H$\index{H@${}_\mathcal{Q}\H=\mathcal{Q}\otimes_\Ze\H$}
denote the
$\mathcal{Q}$-algebra ${}_\mathcal{Q}\H=\mathcal{Q}\otimes_\Ze\H$.
Inside ${}_\mathcal{Q}\H$ we normalize the elements $\i_s$
as follows.

We first introduce
\begin{equation}
n_\a:=(q_{\alpha^\vee}^{1/2}+\theta_{-\alpha/2})
(q_{\alpha^\vee}^{1/2}q_{2\alpha^\vee}-\theta_{-\alpha/2})\in
\A.
\end{equation}

Then the normalized intertwiners $\iota_s^0$
($s\in S_0$) are defined by (with $s=s_\a$,
$\a\in R_1$):
\begin{equation}\label{eq:defint}
\i^0_s :=n_\a^{-1}\i_s\in {}_\mathcal{Q}\H.
\end{equation}
\index{0j@$\iota_s^0$, normalized intertwining element of
${}_\mathcal{Q}\H$}
It is known that the normalized elements $\iota_s^0$
satisfy $(\i_s^0)^2=1$.
In particular, $\i_s^0\in{}_\mathcal{Q}\H^\times$, the
group of invertible elements of ${}_\mathcal{Q}\H$.
In fact we have:
\begin{lemma}(\cite{O}, Lemma 4.1)
The map $S_0\ni s\to \i^0_s\in {}_\mathcal{Q}\H^\times$ extends (uniquely)
to a homomorphism $W_0\ni w\to
\i^0_w\in {}_\mathcal{Q}\H^\times$. Moreover,
for all $f\in {}_\mathcal{Q}\A$ we have that
$\i_w^0f\i_{w^{-1}}^0=f^w$.
\end{lemma}
\subsection{Hilbert algebra structure on $\H$}
The anti-linear map $h\to h^*$ defined by $(\sum_wc_wN_w)^*=\sum_w
\overline{c}_{w^{-1}}N_w$\index{*@$h\to h^*$, conjugate
linear anti-involution of $\H$}
is an anti-involution of $\H$. Thus it
gives $\H$ the structure of an involutive algebra.

In the context of involutive algebras we also dispose
of Schur's lemma for topologically irreducible representations
(cf. \cite{Dix}). Thus the topologically irreducible representations
of the involutive algebra $(\H,*)$ are finite dimensional
by Theorem \ref{thm:cent}.

The linear functional
$\tau:\H\to \mathbb{C}$\index{0t@$\tau$, trace functional of $\H$}
given by
$\tau(\sum_wc_wN_w)=c_e$ is a positive trace for the involutive
algebra $(\H,*)$.
The basis $N_w$ of $\H$ is orthonormal with respect to the
pre-Hilbert structure
$(x,y):=\tau(x^*y)$\index{$(\cdot,\cdot)$!inner product on $\H$}
on $\H$. We denote the
Hilbert completion of $\H$ with respect to $(\cdot,\cdot)$
by $L_2(\H)$\index{L@$L_2(\H)$, Hilbert completion of $\H$}.
This is a separable Hilbert space with Hilbert
basis $N_w$ ($w\in W$).

Let $x\in\H$. The operators $\l(x):\H\to\H$ (given by $\l(x)(y)=xy$)
and $\rho(x):\H\to\H$ (given by $\rho(x)(y):=xy$) extend to
$B(L_2(\H))$, the algebra of bounded
operators on $L_2(\H)$. This gives $\H$ the structure of
a Hilbert algebra (cf. \cite{Dix}).

The operator norm completion of $\l(\H)\subset B(L_2(\H))$
is a $C^*$-algebra which we call the reduced $C^*$-algebra
$C^*_r(\H)$\index{C@$\cf$, the reduced $C^*$ algebra of $\H$}
of $\H$.
The natural action of $C^*_r(\H)$ on $L_2(\H)$ via
$\l$ (resp. $\rho$) is called the left regular (resp. right
regular) representation of $C^*_r(\H)$. Since it has only
finite dimensional irreducible representations by the above remark,
$C^*_r(\H)$ is of type I.

The norm
$\Vert x\Vert_o$ of $x\in C^*_r(\H)$\index{$\Vert\cdot\Vert_o$,
operator norm on $\H$}
is by definition equal to the norm of $\l(x)\in B(L_2(\H))$.
Observe that the map $x\to\l(x)N_e$ defines an embedding
$\cf\subset L_2(\H)$.
\subsection{Discrete series representations}
\begin{dfn}
We call an irreducible
representation $(V_\d,\d)$ of $(\H,*)$
a {\it discrete series representation} if $(V_\d,\d)$
is equivalent to a subrepresentation of $(L_2(\H),\l)$.
We denote by
$\Delta=\Delta_{\Ri,q}$\index{0D3@$\Delta=\Delta_{\Ri,q}$,
complete set of representatives of the irreducible discrete
series representations of $\H(\Ri,q)$}
a complete set of representatives of the equivalence classes of the
irreducible discrete series representations of $(\H,*)$.
When $r\in T$ is given,
$\Delta_{W_0r}\subset \Delta$\index{0D31@$\Delta_{W_0r}$,
subset of $\Delta$ consisting of the
representations with central character $W_0r$}
denotes the subset of $\Delta$ consisting
of irreducible discrete series representations with
central character $W_0r$ ($r\in T$).
\end{dfn}
\begin{cor}\label{cor:fin}
(of Theorem \ref{thm:cent})
$\Delta_{W_0r}$ is a finite set.
\end{cor}
There is an important characterization of the discrete series
representations due to Casselman. This characterization has
consequences for the growth behaviour of matrix
coefficients of discrete series representations.
Recall that $T$ denotes the complex algebraic torus of
characters of the lattice $X$. It has polar decomposition
$T=T_{rs}T_u$\index{T3@$T_u=\operatorname{Hom}(X,S^1)$,
compact form of $T$}\index{T1@$T_{rs}=\operatorname{Hom}(X,\R^\times_+)$,
real split form of $T$}
where $T_{rs}$ is the real split form of $T$,
and $T_u$ the compact form. If $t\in T$ we denote
by
$|t|\in T_{rs}$\index{t9@$\vert t\vert \in T_{rs}$,
real split part of $t\in T$}
its real split part.
\begin{thm}\label{thm:casds} (Casselman's criterion for discrete
series representations, cf. \cite{O}, Lemma 2.22).
Let $(V_\d,\d)$ be an irreducible representation of $\H$.
The following are equivalent:
\begin{enumerate}
\item[(i)] $(V_\d,\d)$ is a discrete series representation.
\item[(ii)] All matrix coefficients of $\d$ belong to $\hf$.
\item[(iii)] The character $\chi$ of $\d$ belongs to $\hf$.
\item[(iv)] The weights $t\in T$ of the generalized
$\A$-weight spaces of $V_\d$ satisfy: $|x(t)|<1$, for all
$0\not=x\in X^+$.
\item[(v)] $Z_X=\{0\}$, and there exists an $\e>0$ such that
for all matrix coefficients $m$ of $\d$, there exists a $C>0$
such that the inequality $|m(N_w)|<C\q^{-\e l(x)}$ holds.
\end{enumerate}
\end{thm}
We have the following characterization of the set of central
characters of irreducible discrete series representations.
For the notion of ``residual points'' of $T$ we refer the reader
to Definition \ref{dfn:res}.
\begin{thm}\label{thm:dsres}
(cf. \cite{O}, Lemma 3.31 and Corollary 7.12)
The set $\Delta_{W_0r}$ is nonempty iff $r\in T$ is a residual
point. In particular, $\Delta$ is finite, and empty unless
$Z_X=0$.
\end{thm}
\subsection{The Schwartz algebra; tempered representations}
We define norms $p_n$
($n\in\Z_+=\{0,1,2,\dots\}$)
on $\H$ by
\begin{equation}
p_n(h)=\max_{w\in W}|(N_w,h)|(1+\nc(w))^n,
\end{equation}
and we define the Schwartz completion
$\S$\index{S@$\S$, the Schwartz completion of $\H$}
of $\H$ by
\begin{equation}
\S:=\{x=\sum_w x_w N_w\in\H^*\mid p_n(x)<\infty\ \forall
n\in\Z_+\}
\end{equation}
In (\cite{O}, Theorem 6.5) it was shown that the multiplication
operation of $\H$
is continuous with respect to the family $p_n$ of norms.
The completion $\S$ is a (nuclear, unital) Fr\'echet algebra
(cf. \cite{O}, Definition 6.6).

As an application of (\cite{O}, Theorem 6.1) it is easy to see
that there exist constants $D\in\Z_+$ and $C>0$ such that
$\Vert h\Vert_o\leq Cp_D(h)$ for all $h\in \H$. Thus
\begin{equation}\label{eq:sinl2}
\S\subset\cf\subset L_2(\H)
\end{equation}
The Main Theorem \ref{thm:main}
can be viewed as a structure theorem for this Fr\'echet algebra
via the Fourier transformation.
\begin{dfn}\label{dfn:temp}
The topological dual $\S^\prime$ is called the space of
tempered functionals. A continuous representation of $\S$
is called a tempered representation. By abuse of terminology,
we call a finite dimensional
representation of $\H$ tempered if it extends
continuously to $\S$.
\end{dfn}
In particular, a finite dimensional representation
$(V,\pi)$ of $\H$
is tempered if and only if the matrix
coefficient $h\to\phi(\pi(h)v)$ extends continuously to $\S$
for all $\phi\in V^*$ and $v\in V$.

We will now discuss Casselman's criteria for
temperedness of finite functionals and finite dimensional
representations if $\H$.
\subsection{Casselman's criteria for temperedness}
\subsubsection{Algebraic dual of $\H$}
We identify the algebraic dual
$\H^*$\index{H1@$\H^*$, algebraic dual of $\H$}
of $\H$ with formal linear
combinations $f=\sum_{w\in W} d_wN_w$ via the sesquilinear pairing
$(\cdot,\cdot)$ defined by $(x,y)=\tau(x^*y)$. Thus $f(x)=(f^*,x)$
and $d_w=f(N_{w^{-1}})$. For $x,y\in\H$ and $f\in\H$ we define
$R_x(f)(y)=f(yx)$ and $L_x(f)(y):=f(xy)$ (a right representation
of $\H$). Note that in terms of multiplication of formal series we
have: $R_x(f)=x.f$ and $L_x(f)=f.x$ (sic).
\subsubsection{Finite functionals}
Let $\Ab\subset\H^*$
denote the linear space of {\it finite}
linear functionals on $\H$:
\begin{dfn}\index{A@$\Ab$, space of finite functionals on $\H$}
The space $\Ab$ consists of all the elements $f\in \H^*$
such that the space $\H.f.\H$ is finite dimensional.
\end{dfn}
Since $\H$ is finite over its center $\Ze$, $f$ is finite if and
only if
$\operatorname{dim}(f.\Ze)<\infty$. Let $\A$ denote
the abelian subalgebra of $\H$ spanned by the elements $\theta_x$
with $x\in X$. Since $\Ze\subset \A$ we see that $f\in \Ab$ if and
only if $\operatorname{dim}(\A.f)<\infty$ if and only if
$\operatorname{dim}(f.\A)<\infty$.
\subsection{Exponents of finite functionals}\label{sub:exp}
\begin{dfn}\label{dfn:exp}
We say that $t\in T$ is an exponent of $f\in\Ab$ if
the $X$-module on the finite dimensional space $V=f.\H$ (the space of
{\it left} translates of $f$) defined via $x\to L_{\theta_x}|_V$
contains a (generalized) weight space with weight $t$.
\end{dfn}
\begin{prop}
Let $f\in\Ab$ and let $\epsilon$ denote the set of exponents of $f$.
There
exist unique functions $E_t^f$ ($t\in\epsilon$) on $\H\times X$,
{\it polynomial} in $X$, such that
\begin{equation}\label{eqn:exp}
f(\theta_xh)=\sum_{t\in\epsilon}E_t^f(h,x)t(x)
\end{equation}
\end{prop}
\begin{proof}
{\it{Uniqueness}}:
Suppose that we have a finite set $\epsilon$ of exponents and
for each $t\in\epsilon$ a polynomial function $x\to E_t(x)$
of $X$ such that
\begin{equation*}
\sum_{t\in\epsilon}E_t(x)t(x)\equiv 0.
\end{equation*}
Suppose that there exists a $t\in \epsilon$
such that $x\to E_t(x)$ has positive degree. We apply the
difference operator $\Delta_{t,y}$ ($t\in\epsilon, y\in X$)
defined by
\begin{equation*}
\Delta_{t,y}(f)(x):=t(y)^{-1}f(x+y)-f(x).
\end{equation*}
It is easy to see that for a suitable choice of $y$ this operator
lowers the degree of the coefficient of $t$ by $1$, and leaves the
degrees of the other coefficients invariant. Hence, if we assume
that not all of the coefficients $E_t$ are zero, we obtain a nontrivial
complex linear relation of characters of $X$, after applying a suitable
sequence of operators $\Delta_{s,z}$. This is a contradiction.

{\it{Existence}}: We fix $h\in\H$ and we decompose $f$ according
to generalized $L_X$-eigenspaces in $V$. We may replace $f$ by one
of its constituents, and thus assume that $\epsilon=\{t\}$.
We may replace the action of $X$ by the action
$L^\prime_x=t(x)^{-1}L_x$. Therefore it is enough to consider the
case $t=1$. Let $N$ denote the dimension of $V$. By Engel's
theorem applied to the commuting unipotent elements
$L_{\theta_x}$ acting in $V$, we see that any product of
$N$ or more difference operators of the form
$\Delta_y=L_{\theta_y}-1$ is equal to zero in $V$.
By induction on $N$ this
implies that for any $h$, the function $x\to f(\theta_x h)$ is
a polynomial in $x$ of degree at most $N-1$.
\end{proof}
\begin{cor}
We have $E_t^f(\theta_x h,y)=t(x)E_t^f(h,x+y)$. In particular,
the degree of the polynomial $E_t^f(h,x)$ is uniformly bounded
as a function of $h$.
\end{cor}
\begin{cor}\label{cor:ft}
Put $f_t(h)=E_t^f(h,0)$. Then $f_t$ is the component of $f$
corresponding to the generalized $L_X$-eigenspace with eigenvalue $t$
in $V$. Observe that $f_t(\theta_x h)=t(x)E_t^f(h,x)$, and that $f_t\in
f\cdot\mathcal{A}=L_X(f)\subset V\subset \Ab$.
\end{cor}
\subsection{The space
$\Ab^{temp}$\index{A1@$\mathbb{A}^{temp}$, finite tempered functionals
on $\H$} of tempered finite functionals}
If $f\in\Ab$, we can express the condition $f\in \S^\prime$
(temperedness) or $f\in L_2(\H)$ (square integrability) in terms
of a system of inequalities on the set of exponents $\epsilon$ of
$f$. This is the content of the Casselman conditions for
temperedness (\cite{O}, Lemma 2.20). We will formulate these
results below, adapted to suit the applications we have in mind
(Section \ref{sec:ct}).

We introduce a partial ordering
$\leq_{F_0}$\index{<@$\leq_P$, partial ordering of exponents
relative to $P\in\P$}
on $T_{rs}$ by
\begin{equation}
t_1\leq_{F_0} t_2 \Longleftrightarrow x(t_1)\leq x(t_2)
\mathrm{\ for\ all\ } x\in X^+
\end{equation}
(this is in fact the special case $P=F_0$ of
the ordering $\leq_{P}$ defined in Definition
\ref{dfn:expordpar}).

Let $(V,\pi)$ be finite dimensional
representation of $\H$.
It follows easily from Definition \ref{dfn:exp}
that the union of the sets of exponents of the matrix coefficients
$h\to\phi(\pi(h)v)$ of $\pi$ coincides with the set of weights $t$ of
the generalized $\A$-weight spaces of $V$.
Using (\cite{O}, Lemma 2.20) we get:
\begin{cor}\label{cor:cas}(\cite{O}, Lemma 2.20)
Let $(V,\pi)$ be a finite dimensional representation of $\H$.
The following statements are equivalent:
\begin{enumerate}
\item[(i)] $(V,\pi)$ is tempered.
\item[(ii)] The weights $t$ of
the generalized $\A$-weight spaces of $V$ satisfy
$|t|\leq_{F_0} 1$.
\item[(iii)] The exponents $t$ of the matrix coefficients
$h\to\phi(\pi(h)v)$ of $\pi$ satisfy $|t|\leq_{F_0} 1$.
\end{enumerate}
\end{cor}
Let $f\in\Ab$. The space of matrix coefficients of the
finite dimensional representation $(V_f:=R_\H(f),R)$ is the
space $\H\cdot f\cdot \H$. Hence the union of the sets
of exponents of the matrix coefficients of $V_f$ is equal
to the set of exponents of $f$. Hence we obtain:
\begin{cor}\label{cor:castemp}
(Casselman's condition)
We have $f\in \Ab^{temp}:=
\Ab\cap\mathcal{S}^\prime$\index{A@$\Ab^{temp}$, tempered elements in $\Ab$}
if and only
if the real part $|t|$ of every exponent $t$ of $f$ satisfies
$|t|\leq_{F_0} 1$.
\end{cor}
\begin{dfn}\label{dfn:cusp}
\index{A@$\Ab^{cusp}$, cuspidal elements in $\Ab^{temp}$}
We put $\Ab^{cusp}$ for the subspace of $\Ab^{temp}$ consisting
of those $f$ such that all exponents $t$ of $f$ satisfy
$|t|=\prod_{\a\in F_0}(d_a\otimes\alpha^\vee)$ with $0<d_\a<1$.
\end{dfn}
Then Theorem \ref{thm:casds} implies similarly that:
\begin{cor}\label{cor:casmcds}
$\Ab_2:=\Ab\cap L_2(\H)\not=0$ only if
$Z_X=0$, and in this case,
$\Ab_2=\Ab^{cusp}$.
\end{cor}
\subsection{Formal completion of $\H$ and Lusztig's structure theorem}
\label{sub:compl}
Let $t\in T$, and let $\I_t$ denote the maximal ideal of $\Ze$
associated with the orbit $W_0t$. We denote by
$\bar{\Ze}_{W_0t}$\index{Z92@$\bar\Ze_{W_0t}$, formal completion
of $\Ze$ at $W_0t$}
the $\I_t$-adic completion of $\Ze$. In \cite{Lu} Lusztig considered
the structure of the completion
\begin{equation}\label{eq:compl}
\index{H6@$\bar{\H}_t=\bar\Ze_{W_0t}\otimes_\Ze \H$}
\bar{\H}_t:=\bar\Ze_{W_0t}\otimes_\Ze \H.
\end{equation}
We will use Lusztig's results on the structure of this
formal completion (in a slightly adapted version)
for so called $R_P$-generic points $t\in T$.
\subsubsection{$R_P$-generic points of $T$}\label{subsub:generic}
Let $R_P\subset R_0$ be a parabolic subset of roots, i.e.
$R_P=\R R_P\cap R_0$.
Let us recall the notion of an $R_P$-generic point $t\in T$
(cf. \cite{O}, Definition 4.12). With $t\in T$ we associate
$R_{P(t)}\subset R_0$, the smallest parabolic subset containing
all roots $\a\in R_0$ for which one of the following statements
holds (where $c_\a$ denotes the Macdonald $c$-function, cf.
equation (\ref{eq:defcr1})):
\begin{enumerate}
\item[(i)] $c_\a\not\in\O_t^\times$ (the invertible holomorphic
germs at $t$)
\item[(ii)] $\a(t)=1$
\item[(iii)] $\a(t)=-1$ and $\a\not\in 2X$.
\end{enumerate}
We say that $t_1,t_2\in T$ are equivalent if there exists a
$w\in W_{P(t_1)}:=W(R_{P(t_1)})$ such that $t_2=w(t_1)$. Notice
that in this case $R_{P(t_1)}=R_{P(t_2)}$, so that this is indeed
an equivalence relation.
The equivalence class of $t\in T$ is
equal to the orbit
$\varpi=W_{P(t)}t\subset W_0t$\index{0p1@$\varpi$, equivalence
class inside $W_0t$}.

We define $P(t)$ as the basis of simple roots of
$R_{P(t)}$ inside $R_{0,+}$, and we sometimes use the
notation $P(\varpi)$ instead of $P(t)$.
\begin{dfn}\label{dfn:gen}
We call $t\in T$ an $R_P$-generic point if $wt\in \varpi$ (with
$w\in W_0$) implies that $w\in W_P$.
\end{dfn}
\begin{rem}\label{rem:nongen}
Notice that if $t\in T$ is $R_P$ generic then $R_{P(t)}\subset R_P$
(but not conversely). In particular, the set of $R_P$ nongeneric
points is contained in a finite union of cosets of the finite
collection of codimension $1$ complex subtori $H$ of $T$ such
that $\a(H)=1$ for some $\a\in R_0\backslash R_P$.
\end{rem}
\subsubsection{Lusztig's first reduction Theorem}
Lusztig \cite{Lu} associates idempotents $e_{w\varpi}\in\bar\H_t$
\index{e1@$e_\varpi$, idempotent of $\bar\H_t$ associated with
$\varpi\subset W_0t$}
with the equivalence classes $w\varpi\in W_0t$.
By Lusztig's first reduction Theorem (cf. \cite{Lu})
we know that if $u,v\in W^P$, then
$\iota_u^0 e_\varpi\iota_{v^{-1}}^0$ is a well defined element
of $\bar\H_t$, and that we have the decomposition (compare with
\cite{O}, equation (4.46))
\begin{equation}\label{eq:dec}
\bar\H_t=\bigoplus_{u,v\in W^P}
\iota_u^0e_\varpi\bar\H^P_t\iota_{v^{-1}}^0.
\end{equation}
Moreover, the subspace $\iota_u^0e_\varpi\bar\H^P_t\iota_{v^{-1}}^0$ is
equal to $e_{u\varpi}\bar\H_t e_{v\varpi}$. When $u=v$ then this is
a subalgebra of $\bar\H_t$, and when $u=v=e$ then this subalgebra reduces to
$e_\varpi\bar\H^P_t$, which is isomorphic to $\bar\H^P_t$ via
$x\to e_\varpi x$. Finally, for $x\in\bar\H^P_t$ we have
the formula
$\iota_u^0(e_\varpi x)\iota_{u^{-1}}^0=e_{u\varpi}\psi_u(x)$.

We will use this in the situation that $t\in T$ is
is of the form $t=r_Pt^P$ with $W_Pr_P\subset T_P$ the central
character of a discrete series representation $(V_\d,\d)$, and
$t^P\in T^P$ (this is the case if $W_0t\subset T$ is the central
character of a representation which is induced from $(V_\d,\d)$.
Recall that in this situation $r_P\in T_P$ is an $(R_P, q_P)$-residual
point (Theorem \ref{thm:dsres}).
Therefore, $R_{P(t)}\supset R_P$ (\cite{O}, Proposition 7.3),
and $R_{P(t)}=R_P$ for an open dense subset of $T^P$ (the
complement of a subvariety of codimension $1$ in $T^P$).
Thus $t=r_Pt^P$ is $R_P$-generic iff $R_{P(t)}=R_P$, and
then the equivalence class of $t$ is equal to
$\varpi=W_Pt$.
\subsubsection{Application}\label{subsub:ext}
We will use the above result (\ref{eq:dec}) when analyzing a
finite functional $f\in\mathbb{A}$ or a representation $\pi$
of $\H$ which contains a power $\I^n_t$ of $\I_t$ in its kernel.

We can thus view $f$ (or $\pi$) as a linear function on the
quotient $\H/\I_t^n\H$. Since this quotient is finite dimensional
(by Theorem \ref{thm:cent}), we have
\begin{equation}\label{eq:inn}
\H/\I_t^n\H=\bar{\H}_t/\I_t^n\bar{\H}_t.
\end{equation}
In this way we can view $f$ (resp. $\pi$) as a functional
(resp. representation) of the completion $\bar{\H}_t$.
For example, this applies when $W_0t$ is the central character
of an irreducible representation $\pi$. We can view $\pi$
as a representation of the quotient $\H^t:=\H/\I_t\H$ (the
case $n=1$ of (\ref{eq:inn}))\index{H7@$\H^t=\H/\I_t\H$,
where $\I_t$ is the maximal ideal of $W_0t$ in $\Ze$},
and the matrix coefficients of $\pi$
can be viewed as functionals on $\H^t$.
\section{Fourier Transform}
In this section we briefly review the Fourier transform on
$L_2(\H)$ as formulated in \cite{O}. The spectral data are
organized in terms of the induction functor on the groupoid
of unitary standard induction data $\W_{\Xi_u}$.
Finally we formulate the Main Theorem
\ref{thm:main} and discuss its applications.
\subsection{Induction from standard parabolic
subquotient algebras}\label{sub:stind}
Let $P\subset F_0$ and let $W_P\subset W_0$ be the standard
parabolic subgroup of $W_0$ generated by the simple reflections
$s_\a$\index{s@$s_\a$, reflection in $\a$}
with $\a\in P$.
Let $\H^P\subset \H$ be the subalgebra
$\H^P:=\H(W_P)\cdot\A\subset \H$, and let $\H_P$ denote the
quotient of $\H^P$ by the (two sided) ideal generated by
the central elements $\theta_x-1$ where $x\in X$ is such
that $\langle x,\a^\vee\rangle=0$ for all $\a\in P$.
Then $\H_P$ is again an affine Hecke algebra,
with root datum $\Ri_P=(R_P,X_P,R_P^\vee,Y_P,P)$, where
$X_P=X/P^{\vee,\perp}$ and $Y_P=Y\cap \R P^\vee$, and
root labels $q_P$ that are obtained by restriction from
$R_{\mathrm{nr}}$ to $R_{P,\mathrm{nr}}$.

There exists a parameter family of homomorphisms
$\phi_{t^P}:\H^P\to \H_P$
with $t^P\in T^P\subset T$\index{T@$T^P\subset T$,
algebraic subtorus of $T$ with
          $\operatorname{Lie}(T^P)=P^\perp$},
the subtorus with character
lattice $X^P=X/(X\cap\R P)$, defined by
$\phi_{t^P}(\theta_x T_w)=x(t^P)\theta_{\bar{x}}T_w$, where
$\bar{x}\in X_P$ denotes the canonical image of $x$ in $X_P$.
The kernel of $\phi_{t^P}$ is the
two-sided ideal generated by elements of the form
$x(t^P)^{-1}\theta_x-1$, with $x\in X$ such that
$\langle x,\a^\vee\rangle=0$ for all $\a\in P$.

Let $(V_\d,\d)$ be a discrete series representation of the
subquotient Hecke algebra $\H_P$.
Let $W_Pr_P$ be the central character of $\d$. It is known that
$r_P$ is a residual point of
$T_P$\index{T@$T_P\subset T$, algebraic
subtorus orthogonal to $T^P$, with character lattice $X_P$}
(cf. \cite{O}, Lemma 3.31),
the subtorus of $T$ with
character lattice $X_P$.

Now let $t^P\in T_u^P$, and let $\d_{t^P}$ denote the lift
to $\H^P$ of $\d$ via $\phi_{t^P}$. Then the induced representation
$\pi=\pi(\Ri_P,W_Pr_P,\d,t^P)$\index{0p@$\pi(\Ri_P,W_Pr,\d,t)=
\operatorname{Ind}_{\H^P}^\H(\d_t)$,
parabolically induced representation}
from the representation
$\d_{t^P}$ of $\H^P$ to $\H$ is a unitary, tempered
representation (cf. \cite{O}, Proposition 4.19 and
Proposition 4.20).
\subsubsection{Compact realization of $\pi(R_P,W_Pr_P,\d,t^P)$}
\label{subsub:comp}
Put $\H(W^P)\subset\H$ for the finite
dimensional linear subspace of $\H$ spanned by the elements
$N_w$ with $w\in W^P$.
Then
\begin{equation}
\H\simeq\H(W^P)\otimes\H^P,
\end{equation}
where the isomorphism is realized by the product map.
Therefore we have the isomorphism
\begin{equation}
\index{i@$i(V_\d)=\H(W^P)\otimes V_\d$, underlying vector space of
$\pi(P,\d,t^P)$ in the compact realization}
\H\otimes_{\H^P}V_\d\simeq i(V_\d):=\H(W^P)\otimes V_\d.
\end{equation}
We will use this isomorphism to identify the representation
space of $\pi(P,W_Pr_P,\d,t^P)$
with $i(V_\d)$. This realization
of the induced representation is called the {\it compact
realization}, by analogy with induced representations for
reductive groups.

According to \cite{O}, Proposition 4.19, the representation
$\pi(P,W_Pr_P,\d,t^P)$ is unitary (i.e. a $*$-representation)
with respect to the Hermitian inner product
\begin{equation}\label{eq:in}
\index{<@$\langle\cdot,\cdot\rangle$!inner product on $V_\pi$}
\langle h_1\otimes v_1,h_2\otimes v_2\rangle
=\tau(h_1^*h_2)(v_1,v_2),
\end{equation}
where $(v_1,v_2)$ denotes the inner product on the representation
space $V_\d$ of the discrete series representation $(V_\d,\d)$.

More generally, for $t^P\in T^P$ the Hermitian form
$\langle\cdot,\cdot\rangle$ on $i(V_\d)$ defines a nondegenerate
sesquilinear pairing of $\H$-modules as
follows:
\begin{equation}
\pi(P,W_Pr_P,\d,\bar{t^P}^{-1})\times
\pi(P,W_Pr_P,\d,t^P)\to\mathbb{C}.
\end{equation}
\subsection{Groupoid of standard induction data}
Let $\P$ denote the power set of $F_0$.
Let
$\Xi$\index{0O@${\Xi}$ set of standard induction data
$\xi=(P,\d,t^P)$ with $P\in\P, \d\in\Delta_P$ and $t^P\in T^P$}
(respectively $\Xi_u$\index{0O2@${\Xi_u}$ set
of tempered standard induction data $\xi=(P,\d,t^P)$})
denote the set of all
triples
$\xi=(P,\d,t^P)$\index{0o@$\xi=(P,\d,t^P)$ standard induction datum}
with
$P\in \P$, $\d$ an irreducible discrete series
representation of $\H_P$
(with underlying vector space $V_\d$), and $t^P\in T^P$
(respectively $t^P\in T^P_u$).
We denote the central character of $\d$ by $W_Pr_P$.

Let
$\W$\index{W8@$\mathcal{W}$, groupoid whose set of objects is
$\mathcal{P}$, with morphisms
$\operatorname{Hom}_\mathcal{W}(P,Q)=
\mathcal{W}(P,Q):=K_Q\times W(P,Q)$}
denote the finite groupoid whose set of objects
is $\P$ and such that the set of arrows from
$P$ to $Q$ ($P,Q\in \P$) consists of $K_Q\times W(P,Q)$, where
$K_Q=T_Q\cap T^Q$\index{K@$K_P$, finite abelian group $T_P\cap T^P$}
and $W(P,Q)=\{w\in W_0\mid w(P)=Q\}$\index{W5a@$W(P,Q)=
\{w\in W_0\mid w(P)=Q\}$,
with $P,Q\in\P$}. The
composition of arrows is defined by
$(k_1,w_1)(k_2,w_2)=(k_1w_1(k_2),w_1w_2)$.
This groupoid acts on $\Xi$ as follows.
An element $g=k\times n\in K_Q\times W(P,Q)$ of $\W_\Xi$
defines
an algebra isomorphism
$\psi_g:\H_P\to\H_Q$\index{0x@$\psi_g:\H^P\to\H^Q$, isomorphism for
$g\in K_Q\times W(P,Q)$}
as follows.
An element $n\in W(P,Q)$ defines an
isomorphism from the root datum $(\Ri_P,q_P)$ to $(\Ri_Q,q_Q)$, which
determines an algebra isomorphism $\psi_n$.
On the other hand, if $k\in K_Q$ then
$\psi_k:\H_Q\to\H_Q$ is the automorphism
defined by $\psi_k(\theta_x N_w)=k(x)\theta_x N_w$.
Then $\psi_g$
is defined by the composition of these isomorphisms.
We obtain
a bijection
$\Psi_g:\Delta_{W_Pr_P}\to\Delta_{k^{-1}W_Qn(r_P)}$\index{0X@$\Psi_g:
\Delta_{P}\to\Delta_{Q}$,
bijection induced by $\psi_g$}
(where $\Delta_{W_Pr_P}=\Delta_{P,W_Pr_P}$\index{0D4@$\Delta_P
=\Delta_{\Ri_P,q_P}$,
complete set of representatives for the equivalence
classes of irreducible discrete
series representations of $\H_P\\=\H(\Ri_P,q_P)$}
denotes a complete set of representatives for the
equivalence classes of irreducible discrete series representations
of $\H_P$ with central character $W_Pr_P$) characterized by the
requirement
$\Psi_g(\d)\simeq\d\circ\psi_{g}^{-1}$.
The action of $\W$ on $\Xi$ is defined by:
$g(P,\d,t^P)=(Q,\Psi_g(\d),g(t^P))$, with $g(t^P):=kn(t^P)$.
\begin{dfn}
The fibred product
$\W_\Xi=\W\times_\P\Xi$\index{W8d@$\W_{\Xi}:=\W\times_\P{\Xi}$,
groupoid of standard induction data}\index{W9d@$\W_{\Xi_u}:=
\W\times_\P{\Xi}$,
groupoid of tempered standard induction data}
is called
the groupoid of standard induction data.
The full compact subgroupoid $\W_{\Xi,u}=
\W\times_\P\Xi_u$ is called the groupoid of tempered standard
induction data.
\end{dfn}
\begin{dfn}\label{dfn:genxi}
An element $\xi=(P,\d,t^P)\in\Xi$ is called generic if
$t=r_Pt^P$ is $R_P$-generic
(cf. Definition \ref{dfn:gen}), where $r_P\in T_P$ is such that
$W_Pr_P$ is the central character of $\d$.
\end{dfn}
The groupoid $\W_{\Xi,u}$ was introduced in \cite{O}
(but was denoted by $\W_\Xi$ there) and plays an important
role in the theory of the Fourier transform for $\H$.
It is easy to
see that $\W_\Xi$ is a smooth analytic, \'etale groupoid,
whose set of objects is equal to $\Xi$.
Thus $\W_\Xi$ is
a union of complex algebraic tori, and
therefore we can speak of polynomial and
rational functions on $\Xi$ and on $\W_\Xi$. This
also applies to the full compact subgroupoid $\W_{\Xi,u}$.

Theorem 4.38 of \cite{O} states that there exists an
induction functor
$\pi:\W_{\Xi,u}\to\operatorname{PRep}_{unit,temp}(\H)$,
where the target groupoid is the category of finite dimensional,
unitary, tempered representations of $\H$ in which the morphisms are
given by unitary intertwining isomorphisms modulo the action of scalars.
The image of $\xi=(P,\d,t^P)\in\Xi_u$ is the representation
$\pi(\xi):=\pi(P,W_Pr_P,\d,t^P)$\index{0p@$\pi(\xi):=
\pi(\Ri_P,W_Pr_P,\d,t^P)$, where
$\xi=(P,\d,t^P)$ and $W_Pr_P$ is the central character of $\d$}
of $\H$, in its compact realization,
as was defined in subsection \ref{subsub:comp}.

The intertwining isomorphism
$\pi(g,\xi):i(V_\d)\to i(V_{\Psi_g(\d)})$\index{0p3@$\pi(g,\xi)$,
unitary intertwining operator,
intertwining $\pi(\xi)$ and $\pi(g\xi)$}
associated with
$g=k\times n\in K_Q\times W(P,Q)$ is the operator
$A(g,\Ri_P,W_Pr_P,\d,t^P)$ which was defined in \cite{O}
(equation (4.82)). In order to explain its construction we
need to use Lusztig's theorem on the structure of the formal completion
of $\H$ at the central character of $\pi(\xi)$
(cf. Subsection \ref{sub:compl}).
The central character of $\pi(\xi)$ ($\xi=(P,\d,t^P)$) is
equal to $W_0t$ with $t=r_Pt^P$, where $W_Pr_P$ denotes the
central character of $\d$.
Recall that we can then extend $\pi(\xi)$ to
the formal completion $\bar{\H}_t$ of $\H$ with respect to the
maximal ideal $\I_t$ of $\Ze$ at $W_0t$ (cf. \ref{subsub:ext}).

First we consider the case where $\xi$ is generic
(Definition \ref{dfn:genxi}, Remark \ref{rem:nongen}).
For $w\in W^P$,
$w\not=e$, the idempotent $e_{w\omega}$
(cf. equation \ref{eq:dec}) vanishes on
$1\otimes V_\d\subset i(V_\d)$, where the action is
through $\pi(\xi)$ (extended to the completion).
Therefore we have the natural isomorphisms of vector spaces:
\begin{align}\label{eq:iso}
i(V_\d)
&\simeq \H\otimes_{\H^P} V_\d\\
&\nonumber\simeq \bar{\H}_t\otimes_{e_\varpi\bar{\H}^P_t}V_\d\\
&\nonumber\simeq\bigoplus_{u\in W^P}\iota_u^0e_\varpi
\otimes V_\d,
\end{align}
where ${e_\varpi\bar{\H}^P_t}\simeq \bar{\H}^P_t$ acts on $V_\d$
via $\d_{t^P}$, extended to the formal completion at the central
character $W_Pt$.
We will often suppress the subscript
${e_\varpi\bar{\H}^P_t}$ of $\otimes$.

Let us now define the unitary
standard intertwining operators $\pi(g,\xi)$ in this
case where $\xi$ is generic.
First we {\it choose}
a unitary isomorphism
$\tilde{\d}_g:V_\d\to V_{\Psi_g(\d)}$\index{0dZ@$\tilde\d_g:
V_\d\to V_{\Psi_g(\d)}$,
unitary isomorphism
intertwining the representations $\d\circ\psi_g^{-1}$ and
$\Psi_g(\d)$}
intertwining the representations $\d\circ\psi_g^{-1}$ and
$\Psi_g(\d)$.
Then we define
\begin{align}\label{dfn:A}
\pi(g,\xi):i(V_\d)&\to i(V_{\Psi_g(\d)})\\
\nonumber
h\otimes v&\to
h\iota_{g^{-1}}^0e_{g\varpi}
\otimes_{e_{g\varpi}\H^{g(P)}_{g(t)}}\tilde{\d}_g(v),
\end{align}
where we use the isomorphism of equation (\ref{eq:iso})
to view the right hand side as an element of $i(V_{\Psi_g(\d)})$.
It follows easily that $\pi(g,\xi)$ is an intertwining operator
between $\pi(\xi)$ and $\pi(g\xi)$.

For general $\xi$ we need the following regularity results from
\cite{O}. The matrix elements of $\pi(g,\xi)$ are
meromorphic in $\xi$, with possible poles
at the nongeneric $\xi$.
However, it was shown in \cite{O}, Theorem 4.33, that for
$R_P$-generic $t=r_Pt^P$, $\pi(g,\xi)$
is unitary with respect to the Hilbert space structures
of $i(V_\d)$ and $i(V_{\Psi_g(\d)})$ (which are
independent of $t^P\in T^P_u$,
cf. equation (\ref{eq:in})). Together
with a description of the locus of the possible
singularities
of $\pi(g,\xi)$ (as a rational function on $\Xi_{P,\d}$,
the set of induction data of the form $(P,\d,t^P)$
with $t^P\in T^P$), this implies (according to a simple
argument, cf. \cite{BCD}, Lemma 8)
that $\pi(g,\xi)$ has only
removable singularities in a tubular neighborhood of $\Xi_{P,\d,u}$
(the subset of triples in $\Xi_{P,\d}$ with $t^P\in T^P_u$).
Thus
$\pi(g,\xi)$ has a unique holomorphic extension to a tubular
neighborhood of $\Xi_{P,\d,u}$. This finally clarifies the definition
of $\pi(g,\xi)$ for general $\xi\in \Xi_{P,\d,u}$ (and in fact in
a ``tubular neighborhood'' of this subset of $\Xi_{P,\d}$).

We conclude with the following summary of the above
\begin{thm}\label{thm:fusm}
The induction functor
$\pi:\W_{\Xi,u}\to\operatorname{PRep}_{unit,temp}(\H)$\index{0p@$\pi$,
induction functor on $\W_{\Xi}$}
is
rational and smooth.
\end{thm}
By this we simply mean that on each component $\Xi_{P,\d,u}$
of $\Xi_u$,
the representations $\pi(\xi)$ can be realized by
smooth rational matrices as a function of $\xi\in\Xi_{P,\d,u}$,
and also the matrices of the $\pi(g,\xi)$ are both rational and
smooth in $\xi\in\Xi_{P,\d,u}$.
We note that the matrices $\pi(\xi;h):=\pi(\xi)(h)$
(for $h\in \H$ fixed) are in fact even polynomial, and that the
matrices $\pi(k,\xi)$ (for $k\in K_P$) are constant.
\subsection{Fourier transform on $L_2(\H)$}
Let $V_\xi$ denote the representation space of $\pi(\xi)$, $\xi\in\Xi$.
Thus $V_\xi=i(V_\d)$\index{V@$V_\xi=i(V_\d)$,
underlying vector space of $\pi(\xi)$}
if $\xi=(P,\d,t^P)$, and this vector space
does not depend on the parameter $t^P\in T^P$. We denote by
$\V_\Xi$\index{V3@$\V_{\Xi}$, trivial fiber bundle over $\Xi$}
the trivial fibre bundle over $\Xi$ whose fibre at $\xi$ is
$V_\xi$, thus
\begin{equation}
\V_\Xi:=\cup_{(P,\d)}\Xi_{P,\d}\times i(V_\d)
\end{equation}
where $\Xi_{P,\d}$ denotes the component of $\Xi$ associated
to $P\in\P$, and $(V_\d,\d)\in \Delta_P$, a complete set of
representatives of the irreducible discrete series
representations $\d$ of $\H_P$.
We denote by $\operatorname{End}(\V_\Xi)$ the endomorphism
bundle of $\V_\Xi$, and by
$\operatorname{Pol}(\Xi,\operatorname{End}(\V_\Xi))$\index{Pol@
$\operatorname{Pol}(\Xi)$, space of Laurent polynomials
on $\Xi$}\index{Pol@$\operatorname{Pol}
(\operatorname{End}(\V_\Xi))^\W$,
space of $\W_\Xi$-equivariant sections in
$\operatorname{Pol}(\operatorname{End}(\V_\Xi))$}
the space of
polynomial sections in this bundle. Similarly, let us introduce
the space
$\operatorname{Rat}^{reg}(\Xi_u,\operatorname{End}
(\V_\Xi))$\index{RatO@$\operatorname{Rat}^{reg}(\Xi)$,
rational functions on $\Xi$, regular in an
open neighborhood of $\Xi_u\subset \Xi$}
of rational sections which are regular in a neighborhood of
$\Xi_u$.

There is an action of $\W$ on $\operatorname{End}(\V_\Xi)$ as
follows. If
$(P,g)\in\W_P$\index{W991@$\W_P$, set of elements of $\W$
with source $P\in\P$}
(the set of elements of $\W$ with
source $P\in\P$)
with $g=k\times n\in K_Q\times W(P,Q)$,
$\xi\in\Xi_P$\index{0O1@$\Xi_P$, set of standard induction data
of the form $(P,\d,t^P)$ with $\d\in\Delta_P$ and $t^P\in T^P$},
and $A\in \operatorname{End}(V_\xi)$ we define
$g(A):=\pi(g,\xi)\circ A\circ \pi(g,\xi)^{-1}\in
\operatorname{End}(V_{g(\xi)})$. A section of
$f$ of $\operatorname{End}(\V_\Xi))$ is called
$\W$-equivariant if we have $f(\xi)=g^{-1}(f(g(\xi)))$
for all $\xi\in \Xi$ and $g\in\W_\xi$ (where
$\W_\xi:=\W_P$\index{W992@$\W_\xi=\W_P$ if $\xi\in\Xi_P$}
if $\xi=(P,\d,t^P)$).
\begin{dfn}
We define an averaging projection
$p_\W$\index{p99@$p_\W$, averaging operator for the action of
$\W$ on sections of $\operatorname{End}(\V_\Xi)$}
onto the space of $\W$-equivariant sections by:
\begin{equation}
p_\W(f)(\xi):=|\W_\xi|^{-1}
\sum_{g\in \W_\xi}g^{-1}(f(g(\xi))).
\end{equation}
Notice that this projection preserves the space
$\operatorname{Rat}^{reg}(\Xi_u,\operatorname{End}(\V_\Xi))$,
but not the space
$\operatorname{Pol}(\Xi,\operatorname{End}(\V_\Xi))$.
\end{dfn}
The Fourier transform
$\F_\H$\index{F@$\F$, Fourier Transform!a@$\F_\H$, $\F$
restricted to $\H$}
on $\H$ is the following algebra
homomorphism
\begin{align}\label{eq:ft}
\F_\H:\H&\to
\operatorname{Pol}(\Xi_u,\operatorname{End}(\V_\Xi))^{\W}\\
\nonumber h&\to\{\xi\to\pi(\xi;h)\}
\end{align}
where $\operatorname{Pol}(\Xi,\operatorname{End}(\V_\Xi))^\W$
denotes the space of $\W$-equivariant polynomial sections
of $\operatorname{End}(\V_\Xi))$.

We will now describe a $\W$-invariant measure $\mu_{Pl}$ on
$\Xi_u$ whose push forward to $\W\backslash\Xi_u$ will be the
Plancherel measure of $\H$ (\cite{O}, Theorem 4.43).
Put $\xi=(P,\d,t^P)\in\Xi_u$ and let $t=r_Pt^P$. We write
$d\xi:=|K_{P,\d}|dt^P$ where $dt^P$ denotes the normalized Haar measure
of $T^P_u$ and where $K_{P,\d}$ denotes the stabilizer of $\d$ under
the natural action of $K_P$ on $\Delta_P$. Let
$\K\lhd \W$\index{K@$\K\lhd \W$,
normal subgroupoid whose set of objects is $\P$, with
$\operatorname{Hom}_\K(P,Q)\\=K_P$ if $P=Q$ and $\emptyset$ else}
denote the
normal subgroupoid whose set of objects is $\P$, and with
$\operatorname{Hom}_\K(P,Q)=\emptyset$ unless $P=Q$, in which case we have
$\operatorname{Hom}_\K(P,P)=K_P$. Thus $\W_P/\K_P=\{w\in W_0\mid w(P)\subset
F_0\}$. Let $\mu_{\Ri_P,Pl}(\{\d\})$ denote the Plancherel
mass of $\d$ with respect to $\H_P$ (see \cite{O}, Corollary
3.32 for a product formula for $\mu_{\Ri_P,Pl}(\{\d\})$).
We now define the Plancherel measure $\mu_{Pl}$:
\begin{dfn}\label{dfn:plameas}
\begin{equation}
\index{0m@$\mu_{Pl}$, Plancherel measure, spectral measure
of the tracial state $\tau$ on $\cf$}
d\mu_{Pl}(\xi):=q(w^P)^{-1}|\W_P/\K_P|^{-1}
\mu_{\Ri_P,Pl}(\{\d\})|c(\xi)|^{-2}d\xi
\end{equation}
where $c(\xi)$ is the Macdonald $c$-function, see
Definition \ref{dfn:cind}.
\end{dfn}
This measure is smooth on $\Xi_u$
(Proposition \ref{lem:csmo}(v)),
and it is invariant for the action of $\W$ on $\Xi_u$,
by Proposition \ref{lem:csmo}(ii).

With these notations we have:
\begin{thm}(\cite{O}, Theorem 4.43)\label{thm:mainO}
\begin{enumerate}
\item[(i)]
$\F_\H$ extends to an isometric isomorphism
\begin{equation}\index{F@$\F$, Fourier Transform}
\F:L_2(\H)\to L_2(\Xi_u,\operatorname{End}(\V_\Xi),\mu_{Pl})^\W,
\end{equation}
where the Hermitian inner product $(\cdot,\cdot)$ on
$L_2(\Xi_u,\operatorname{End}(\V_\Xi),
\mu_{Pl})^\W$\index{$(\cdot,\cdot)$!inner product on
$\operatorname{Pol}(\operatorname{End}(\V_{(\K\backslash\Xi)}))$}
is defined by integrating the Hilbert-Schmidt form
$(A,B):=\operatorname{tr}(A^*B)$ in the fibres
$\operatorname{End}(V_\xi)$ against the above measure
$\mu_{Pl}$ on the base space $\Xi_u$.
\item[(ii)] If $x\in\cf\subset L_2(\H)$ then
$\F(x)\in C(\Xi_u,\operatorname{End}(\V_\Xi))^\W$.
\item[(iii)]
Let $\cf^o$ denote the opposite $C^*$-algebra of $\cf$.
Let $(x,y)\in\cf\times\cf^o$ act on $L_2(\H)$ via the regular
representation $\l(x)\times\rho(y)$, and on
$L_2(\Xi_u,\operatorname{End}(\V_\Xi),\mu_{Pl})^\W$ through
fibrewise multiplication from the left with $\F(x)$
and from the right with $\F(y)$. Then $\F$ intertwines
these representations of $\cf\times \cf^o$.
\end{enumerate}
\end{thm}
\begin{proof} As to (ii), first recall that
according to Equation (\ref{eq:ft}), $\F_\H(\H)\subset
\operatorname{Pol}(\Xi_u,\operatorname{End}(\V_\Xi))^{\W}$.
By (\cite{O}, Theorem 4.43(iii)) one easily deduces that
$\Vert h\Vert_o=\Vert\F_\H(h)\Vert_{\operatorname{sup}}$
for all $h\in\H$, where
$\Vert\sigma\Vert_{\operatorname{sup}}:=
\operatorname{sup}_{\xi\in\Xi_u}\Vert\sigma(\xi)\Vert_o$
(where $\Vert\sigma(\xi)\Vert_o$ denotes the operatornorm
of $\sigma(\xi)\in\operatorname{End}(V_\xi)$).
Hence $\F(\cf)\subset C(\Xi_u,\operatorname{End}(\V_\Xi))^\W$.

Now (iii) follows from (ii) and (\cite{O}, Theorem 4.43(iii)).
\end{proof}
The following easy corollary is important in the sequel:
\begin{cor}(\cite{O}, Theorem 4.45)\label{cor:wave}
The averaging operator $p_\W$ defines an orthogonal projection
onto the space of $\W$-equivariant sections in
$L_2(\Xi_u,\operatorname{End}(\V_\Xi),\mu_{Pl})$.
Moreover, if
\begin{equation}\index{J@$\J$, wave packet operator,
adjoint of $\F$}
\J: L_2(\Xi_u,\operatorname{End}(\V_\Xi),\mu_{Pl})\to L_2(\H)
\end{equation}
denotes the adjoint of $\F$ (the wave packet operator),
then $\J\F=\operatorname{id}$ and $\F\J=p_\W$.
\end{cor}
\begin{proof}
Theorem \ref{thm:mainO} implies that $\J\F:=\operatorname{id}$
and that $\F\J$ is equal to the orthogonal projection onto the space of
$\W$-equivariant $L_2$-sections of $\operatorname{End}(\V_\Xi)$.

On the other hand,
since the action of $\W$ on $\operatorname{End}(\V_\Xi)$ is defined
in terms of invertible smooth matrices (cf. Theorem
\ref{thm:fusm}), $p_\W$ preserves the space of $L_2$-sections.
By the $\W$-invariance of $\mu_{Pl}$, the projection $p_\W$ on
$L_2(\Xi_u,\operatorname{End}(\V_\Xi),\mu_{Pl})$ is in fact an
orthogonal projection. This finishes the proof.
\end{proof}
\section{Main Theorem and its applications}\label{sub:mainapp}
The space of smooth section of the trivial bundle
$\operatorname{End}(\V_{\Xi})$ on $\Xi_u$ will be denoted by
$\cc(\Xi_u,\operatorname{End}(\V_{\Xi}))$\index{C@$\cc(\Xi_u,
\operatorname{End}(\V_{\Xi}))$,
space of smooth sections in $\operatorname{End}(\V_{\Xi})$}.
We equip this vector space with its usual Fr\'echet topology.
The collection of semi-norms inducing the topology is of the form
$p(\sigma):=
\operatorname{sup}_{\xi\in\Xi_u}\Vert D\sigma(\xi)\Vert_o$,
where $D$ is a constant coefficient differential
operator on $\Xi_u$ (i.e. one such operator for each
connected component of $\Xi_u$), acting entrywise on the section
$\sigma$
of the trivial bundle $\operatorname{End}(\V_\Xi)$, and where
$\Vert\cdot\Vert_o$ denotes the operatornorm.
It is obvious from the product rule for differentiation that
$\cc(\Xi_u,\operatorname{End}(\V_{\Xi}))$ is a Fr\'echet algebra.

The projection $p_\W$ is continuous on
$\cc(\Xi_u,\operatorname{End}(\V_{\Xi}))$, since it is defined
in terms of the action of $\W$ on $\Xi_u$, and conjugations
with invertible smooth matrices.
Thus the subalgebra $\cc(\Xi_u,\operatorname{End}(\V_{\Xi}))^\W$ of
$\W$-equivariant sections is a closed subalgebra.

We now define the vector space
\begin{dfn}\label{dfn:calC}\index{C@$\mathcal{C}
(\Xi_u,\operatorname{End}(V_{\Xi}))
=c\cc(\Xi_u,\operatorname{End}(V_{\Xi}))$}
\begin{equation}
\mathcal{C}(\Xi_u,\operatorname{End}(V_{\Xi})):=c
\cc(\Xi_u,\operatorname{End}(\V_{\Xi})),
\end{equation}
where $c$ denotes the $c$-function of Definition \ref{dfn:cind}
on $\Xi_u$. We equip
$\mathcal{C}(\Xi_u,\operatorname{End}(V_{\Xi}))$ with the
Fr\'echet space topology of
$\cc(\Xi_u,\operatorname{End}(V_{\Xi}))$ via the
linear isomorphism
$\cc(\Xi_u,\operatorname{End}(V_{\Xi}))\to
\mathcal{C}(\Xi_u,\operatorname{End}(V_{\Xi}))$
defined by $\sigma\to c\sigma$.
\end{dfn}
\begin{lem}\label{lem:ccc}
The complex vector space
$\mathcal{C}(\Xi_u,\operatorname{End}(V_{\Xi}))$
is closed for taking (fibrewise) adjoints, and
\begin{equation}
\mathcal{C}(\Xi_u,\operatorname{End}(V_{\Xi}))
\subset L_2(\Xi_u,\operatorname{End}(V_{\Xi}),\mu_{Pl}).
\end{equation}
Moreover,
\begin{equation}
\cc(\Xi_u,\operatorname{End}(V_{\Xi}))
\subset\mathcal{C}(\Xi_u,\operatorname{End}(V_{\Xi}))
\end{equation}
is a closed subspace.
\end{lem}
\begin{proof}
It is closed for taking adjoints by Proposition
\ref{lem:csmo}(iv) (applied to $d=w^P\in W(P,P^\prime)$),
and it is a subspace of
$L_2(\Xi_u,\operatorname{End}(V_{\Xi}),\mu_{Pl})$
by Proposition \ref{lem:csmo}(i).
The last assertion follows from Proposition
\ref{lem:csmo}(v).
\end{proof}

Now we are prepared to formulate the main theorem of this paper.
\begin{thm}\label{thm:main}
The Fourier transform restricts to an isomorphism of
Fr\'echet algebras
\begin{equation}\label{eq:fs}\index{F@$\F$, Fourier Transform!b@$\F_\S$,
$\F$ restriced to $\S$}
\F_\S:\S\to C^\infty(\Xi_u,\operatorname{End}(\V_\Xi))^\W.
\end{equation}
The wave packet operator $\J$ restricts to a surjective
continuous map
\begin{equation}\label{eq:js}
\J_{\mathcal{C}}:\mathcal{C}(\Xi_u,\operatorname{End}(\V_\Xi))\to\S.
\end{equation}\index{J@$\J$, wave packet operator,
adjoint of $\F$!$\J_{\mathcal{C}}$, restriction of
$\J$ to the Fr\'echet space
$\mathcal{C}(\Xi_u,\operatorname{End}(\V_\Xi))$}
We have $\J_{\mathcal{C}}\F_\S=\operatorname{id_\S}$, and
$\F_\S\J_{\mathcal{C}}=p_{\W,\mathcal{C}}$ (the restriction of
$p_\W$ to $\mathcal{C}(\Xi_u,\operatorname{End}(\V_\Xi))$).
In particular, the map $p_{\W,\mathcal{C}}$ is a continuous
projection of
$\mathcal{C}(\Xi_u,\operatorname{End}(V_{\Xi}))$ onto
$\cc(\Xi_u,\operatorname{End}(V_{\Xi}))^\W$.
\end{thm}
\subsection{Applications of the Main Theorem}
Before we embark on its proof we discuss some
immediate consequences of the Main Theorem.
\begin{cor}(Harish-Chandra's completeness Theorem, cf. \cite{HC2},
and \cite{K}, Theorem 14.31)\label{cor:a}
Let $\xi\in\Xi_u$. The complex linear
span $C_\xi$ of the set of operators \{$\pi(g,\xi)\mid g\in
\operatorname{End}_{\W_\Xi}(V_\xi)\}$ is a unital, involutive
subalgebra of $\operatorname{End}(V_\xi)$.
For all $\xi\in\Xi_{u}$, the centralizer
algebra $\pi(\xi,\H)^\prime$ is equal to $C_\xi$.
\end{cor}
\begin{proof} Let $\xi=(P,\d,t^P)$ and
denote by $C_\xi\subset\operatorname{End}(i(V_\d))$ the complex
linear span of the set of operators
$\{\pi(g,\xi)\mid g\in\operatorname{End}_{\W_\Xi}(V_\xi)\}$.
By Theorem \ref{thm:fusm}, $C_\xi$ is an involutive
(i.e. $*$-invariant), unital subalgebra of
$\operatorname{End}(i(V_\d))$.
Theorem \ref{thm:main} implies that
$\pi(\xi,\H)=C^\prime_\xi$.
The Bicommutant Theorem therefore implies that
$C_\xi=\pi(\xi,\H)^\prime$.
\end{proof}
\begin{cor}\label{cor:b}
The center $\Ze_\S$ of $\S$ is, via the Fourier Transform $\F_\S$,
isomorphic to the algebra $C^\infty(\Xi_u)^\W$.
\end{cor}
\begin{proof}
The algebra of scalar sections of
$C^\infty(\Xi_u,\operatorname{End}(\V_{\Xi_u}))^\W$ is
isomorphic to $C^\infty(\Xi_u)^\W$, and is contained in
$\F_\S(\Ze_\S)$ by Theorem \ref{thm:main}. To show the equality,
observe that Corollary \ref{cor:a} implies that an
element of $\F_\S(\Ze_\S)$ is scalar at all
fibers $\operatorname{End}(V_\xi)$ with $\xi\in\Xi_u$
generic (since $\operatorname{End}_{\W_\Xi}(V_\xi)=\C$ in this case).
By the density of the set of generic points in $\Xi_u$
we obtain the desired equality.
\end{proof}
Notice that $\Ze_\S$ is in general larger than the closure in
$\S$ of the center $\Ze$ of $\H$.
\begin{cor}(Langlands' disjointness Theorem,
cf. \cite{K}, Theorem 14.90)\label{cor:c}
Let $\xi,\xi^\prime\in\Xi_u$. If $\pi(\xi)$ and $\pi(\xi^\prime)$
are not disjoint, then the objects $\xi$, $\xi^\prime$ of
$\W_{\Xi_u}$ are isomorphic (and thus, $\pi(\xi)$ and
$\pi(\xi^\prime)$ are actually equivalent).
\end{cor}
\begin{proof}
Corollary \ref{cor:b} implies that $\Ze_\S$ separates the
$\W$-orbits of $\Xi_u$. Whence the result.
\end{proof}


\begin{cor}\label{cor:e}
The Fourier Transform $\F$ restricts to a $C^*$-algebra
isomorphism
\begin{equation}\index{F@$\F$, Fourier Transform!c@$\F_C$,
$\F$ restricted to $\cf$}
\F_C:\cf\to C(\Xi_u,\operatorname{End}(\V_\Xi))^\W,
\end{equation}
where $\cf$ denotes the reduced $C^*$-algebra of
$\H$ (cf. \cite{O}, Definition 2.4).
\end{cor}
\begin{proof}
By Theorem \ref{thm:mainO}, the restriction of $\F$ to $\cf$
is an algebra homomorphism. It is a homomorphism
of involutive algebras since $\pi(\xi;x^*)=\pi(\xi;x)^*$
(cf. Subsection \ref{sub:stind}).

The reduced $C^*$-algebra $\cf$ of $\H$ is defined in \cite{O}
as the norm closure of $\l(\H)\subset B(L_2(\H))$. By Theorem
\ref{thm:mainO}, the norm
$\Vert x\Vert_o$\index{$\Vert\cdot\Vert_o$, operator norm on $\H$}
of $\cf$ is equal to the
supremum norm $\Vert\F(x)\Vert_{\operatorname{sup}}$ of the
$\W$-invariant continuous function
$\xi\to\Vert\pi(\xi;x)\Vert_o$ on $\Xi_u$
(where $\Vert\pi(\xi;h)\Vert_o$
denotes the operator norm for operators on the finite dimensional
Hilbert space $V_\xi=i(V_\d)$).
Notice that, by the regularity of the standard intertwining operators,
the projection operator $p_\W$ restricts to a continuous
projection on the space of continuous sections of
$\operatorname{End}(\V_\Xi)$.

By Theorem \ref{thm:main}, the closure of $\F(\S)$
with respect to $\Vert\cdot\Vert_{\operatorname{sup}}$ is equal to
$C(\Xi_u,\operatorname{End}(\V_\Xi))^\W$. In view of Theorem
\ref{thm:mainO}(ii) this finishes the proof.
\end{proof}
\begin{cor}
The set of minimal central idempotents of $\cf$ is parameterized
by the (finite) set of $\W$-orbits of pairs $(P,\d)$ with
$P\in\P$ and $\d\in\Delta_P$. The central idempotents $e_{(P,\d)}$
are elements of $\S$.
\end{cor}
\begin{proof}
This is immediate from Theorem \ref{thm:main} and
Corollary \ref{cor:e}.
\end{proof}
\begin{cor}\label{cor:h}
The dense subalgebra $\S\subset\cf$ is closed for holomorphic
functional calculus.
\end{cor}
\begin{proof}
According to a well known criterium for closedness under
holomorphic functional calculus,
we need to check that if $a\in\S$ is invertible in $\cf$,
then $a^{-1}\in\S$. This is obvious from Theorem \ref{thm:main}
and \ref{cor:e}.
\end{proof}
\section{Constant terms of matrix coefficients of
$\pi(\xi)$}
In the remainder of this paper we will prove the Main Theorem
Theorem \ref{thm:main}. A main tool is the notion
of the constant term $f^P$ of a function $f\in\Ab^{temp}$
with respect to a standard parabolic subset $P\in \P$.
\subsection{Definition of the constant terms of
$f\in\Ab^{temp}$}\label{sec:ct}
It this section we define the constant term of
a tempered finite functional
$f\in\Ab^{temp}$ along a standard parabolic subalgebra
$\H^P$ of $\H$. Recall the notion of exponents \ref{sub:exp} and
the Casselman criteria for tempered finite functionals.
\begin{dfn}\label{dfn:expordpar}
\index{<@$\leq_P$, partial ordering of exponents
relative to $P\in\P$}
Let $P\subset F_0$, and $R_P$ the standard parabolic
subsystem with that subset. For {\it real} characters $t_1$, $t_2$
on $X$ we say that $t_1\leq_P t_2$ if and only if $t_1(x)\leq
t_2(x)$ for all $x\in X^{P,+}:=\{x\in X \mid \forall \a\in P:
\langle x,\alpha^\vee\rangle\geq 0\}$.
In other words, $t_1\leq_P t_2$ iff both
$t_1\leq_{F_0} t_2$ and $t_1|_{X\cap P^\perp}=t_2|_{X\cap
P^\perp}$.
\end{dfn}
Thus $t_1\leq_P t_2$ iff $t_1t_2^{-1}=\prod_{\a\in
P}(d_\a\otimes\alpha^\vee)$ with $0<d_\a\leq 1$, where $d\otimes
\a^\vee\in T_{rs}$ is the real character defined by $d\otimes
\a^\vee(x)=d^{\langle x,\a^\vee\rangle}$.

\begin{dfn}(Constant term)\label{dfn:ct}
Let $P\in \P$ and $f\in\Ab^{temp}$. Then we define the
constant term of $f$ along $P$ by
\begin{equation*}
f^P(h):=\sum_{t\in\epsilon:|t|\leq_P 1}f_t(h).
\end{equation*}
where (in the notation of  Corollary \ref{cor:ft})
$f_t(h):=E_t^f(h,0)$, and
the coefficients $E_t^f$ and the set $\epsilon$
are defined by the expansion
\ref{eqn:exp}. We say that an exponent $t\in\epsilon$ of $f$
is $P$-tempered if it satisfies the condition $|t|\leq_P 1$.
\end{dfn}
The notion of cuspidality (cf. Definition \ref{dfn:cusp})
can now be reformulated as follows:
\begin{cor}\label{cor:cuspreform}
Let $f\in\Ab^{temp}$. Then $f\in\Ab^{cusp}$ iff $f^P=0$
for every proper $P\in\P$.
\end{cor}
Observe the following elementary properties of the constant term:
\begin{cor}\label{cor:basic}
\begin{enumerate}
\item[(i)] $f^P\in\Ab^{temp}$.
\item[(ii)] $L_x$ commutes with $f\to f^P$ if $x\in\H^P$.
\item[(iii)] $R_y$
commutes with $f\to f^P$ for all $y\in\H$.
\item[(iv)] $f^P\in L_X(f)=f\cdot\A\subset f\cdot \H$.
\end{enumerate}
\end{cor}
The projection of $f$ to $f^P$ can be made explicit using
an idempotent $e^P$ in a formal completion of $\A\subset\H$.
Such completions were introduced and studied by Lusztig \cite{Lu}
(cf. Subsection \ref{sub:compl}).
This will be applied to the case were $f$ is a matrix coefficient
of a parabolically induced representation in the next subsection.
\subsection{Constant terms of coefficients of $\pi(\xi)$ for
$\xi\in\Xi_u$ generic}
{\it In this subsection we assume that $\xi$ is generic unless
stated otherwise}

We will now discuss the constant terms of a matrix coefficient
of $\pi=\pi(\xi)$ in the case where $\xi=(P,\d,t^P)\in \Xi_u$
is generic.
Choose $r_P\in T_P$ such that $W_Pr_P$ is the central character
of $\d$. We thus assume that $t=r_Pt^P\in T$ is $R_P$-generic in
this subsection.

Let $a,b\in i(V_\d)$, and
denote by $f=f_{a,b}=f_{a,b}(\xi)$ the matrix coefficient
defined by $f(h)=\langle a,\pi(\xi;h)(b)\rangle$.
By \cite{O}, Lemma 2.20 and Proposition 4.20, we have: If
$t^P\in T^P_u$ then
$f_{a,b}\in\Ab^{temp}$ for all $a,b\in
\H(W^P)\otimes V_\d$.
More precisely:
\begin{prop}
The exponents of $f$ are of the form $wt^\prime$ where $w$ runs
over the set $W^P$ and where $t^\prime$ runs over the set of
weights of $\delta_{t^P}$, thus $t^P$ times the set of
$X_P$-weights of $\d$.
\end{prop}
Now let
$Q\subset F_0$ be another standard parabolic. By proof of \cite{O},
Proposition 4.20 we deduce:
\begin{prop}
Let $w\in W^P$ and let $u\in W_P$ such that $wut$ is an exponent
of $f$. If $wu|t|\leq_Q 1$, then $w(P)\subset R_{Q,+}$.
\end{prop}
\begin{proof}
The equivalence class $\varpi$ of $t=r_Pt^P$ is equal to $W_Pt$
(since we assume genericity).
If $wut$ is an exponent
then $wut=w^\prime t^\prime$ with $t^\prime$ an $X$-weight of
$\d_{t^P}$ and $w^\prime\in W^P$.
Thus $t^\prime$ and $ut$ are both in
$\varpi$, the equivalence class of $t$. Hence by genericity,
$w^\prime=w$ and thus $ut=t^\prime$, a weight of $\delta_{t^P}$.
But $\d$ is discrete series for $\H_P$, hence $|ut|=\prod_{\a\in P}
d_\a\otimes\a^\vee$ with all $0<d_\a<1$. Thus $wu|t|=\prod_{\a\in
P}d_\a\otimes w(\a^\vee)\leq_Q 1$ implies (since for all $\a\in
P$: $w(\a^\vee)\in R_{0,+}^\vee$) that $w(P)\subset R_{Q,+}$.
\end{proof}
\begin{cor}
Recall that the equivalence classes in $W_0t$ are of the form
$w\varpi$ with $\varpi=W_Pt$ and $w\in W^P$. If an exponent
$wt^\prime$ of $f$ (with $w\in W^P$ and $t^\prime$ a weight of
$\d_{t^P}$) is $Q$-tempered, then all exponents
of $f$ in its class $w\varpi$ are $Q$-tempered. The class
$w\varpi$ ($w\in W^P$) is $Q$-tempered if and only if
$w(P)\subset R_{Q,+}$.
\end{cor}
\begin{proof}
Since $w(P)\subset R_{Q,+}$ we have $wW_Pw^{-1}\subset W_Q$. Hence
$w\varpi\subset W_Qwt$, so that all elements of $w\varpi$ have the
same restriction $1$ to $X\cap Q^\perp$.
\end{proof}

Now we will express the constant term of
a matrix coefficient of $\pi(\xi)$ in terms of the idempotents
$e_{\varpi}$ of the completion $\bar{\H}_t$. Recall the material
of Subsection \ref{sub:compl}.

We will use the analog of Lusztig's  first reduction
theorem (\ref{eq:dec}) for $\bar{\H}_t$,
in combination with the results in (\cite{O}, Section 4.3)
on the Hilbert algebra structure of
$\overline{\H^t}$\index{H7@$\overline{\H^t}$, the quotient
of $\H^t$ by the
radical of the positive semi-definite Hermitian pairing
$(x,y)_t:=\chi_t(x^*y)$},
the quotient
of $\H^t$ by the radical of the positive semi-definite
Hermitian pairing $(x,y)_t:=\chi_t(x^*y)$,
in order to express and study the constant term.
\begin{prop}\label{cor:idmp}
We have that
\[
f^Q(h)=\sum_{w\in W^P:w(P)\subset R_{Q,+}} f(e_{w\varpi} h)
\]
\end{prop}
\begin{proof}
Let us denote by $J_w$ the ideal in $\A^{W_{w(P)}}$ of elements in
this ring vanishing at $w\varpi$. Clearly $\I_t\subset J_w$ for all
$w$. By some elementary algebra (similar to proof of Prop 2.24(4)
in \cite{EO}) we see
that for every $x\in J_w$ and $k\in\N$ there exist a $\bar{x}\in \I_t$
and a unit $e$ in $\bar{\A}_{m_w}$ such that
\[
ex\in \bar{x}+m_w^k,
\]
where $m_w$ denotes the ideal of all functions in $\A$ which
vanish at the points of $w \varpi$. (To be sure, we construct $\bar{x}$
by first adding an element $u\in J_w^k$ such that $x+u$ is nonzero
at the other classes $w^\prime\varpi$ with $w^\prime\in W^P$,
$w^\prime\not=w$. Take $\bar{x}$ equal to the product of the translates
$(x+u)^w$ where $w$ runs over the set of left cosets $W_0/W_{w(P)}$.
Let $e$ be equal to the product of these factors $(x+u)^w$ where
$w$ runs over the set of left cosets $W_0/W_{w(P)}$ with
$w\not=W_{w(P)}$).
Let $M$ be the ideal of functions in
$\A$ vanishing at $W_0t$. Then $M=\prod m_w$ and by genericity
the ideals $m_w$ are relatively prime. So
$\bar{\A}_M=\oplus\bar{\A}_{m_w}$ by the Chinese remainder
theorem. Then $e_{w\varpi}$ is the unit of the summand
$\bar{\A}_{m_w}$. Let $\bar{e}_{w\varpi}$ be the unit of
$R:=\bar{\A}_{m_w}/\I_t\bar{\A}_{m_w}$ (the canonical image of
$e_{w\varpi}\in \bar{\A}_{m_w}$).
Note that $R$ is finite dimensional over $\mathbb{C}$, and thus
$R$ is Artinian. By definition of $m_w$, $m_w\bar{e}_{w\varpi}$
is contained in all the maximal ideals of $R$. Hence
$m_w\bar{e}_{w\varpi}$ is contained in the intersection of
the maximal ideals of $R$, which is nilpotent in $R$ (see
the proof of Theorem 3.2 of \cite{Mat}). In particular,
for sufficiently large $k$,
$m_w^ke_{w\varpi}\subset \I_t\bar{\A}_{m_w}$, whence
\[
xe_{w\varpi}\in \I_t\bar\A_{m_w}.
\]
But then the right hand side is in the kernel of $\pi$, thus we
conclude that $J_we_{w\varpi}$ is in the kernel of $\pi$. In
particular, the element (for any $z\in X$) $\Theta_z:=\prod_{y\in
W_{w(P)}z} ((wt)(z)^{-1}\theta_y-1)\in J_w$ acts by zero on the
finite dimensional space of left $\A$-translates of $h\to
f(e_{w\varpi} h)$. Thus the exponents of $f\to f(e_{w\varpi} h)$
are contained in $w\varpi$.

We obviously have
\[
f(h)=\sum_{w\in W^P} f(e_{w\varpi}h)
\]
(splitting of $1$ according the decomposition of $\bar{\A}_{M}$).
By the results in this paragraph, an exponent of $h\to
f(e_{w\varpi} h)$ is $Q$-tempered iff all exponents of this
term are $Q$-tempered iff $w(P)\subset R_{Q,+}$. Hence the result.
\end{proof}
\begin{cor}\label{cor:rat}
The constituents $f(e_{w\varpi}h)$ depend on the induction
parameter $t^P$ as a rational function.
\end{cor}
\begin{proof}
In the proof of Corollary \ref{cor:idmp} we can equally well
work over the field $K$
of rational functions on $T^P$ instead of $\mathbb{C}$.
Then $\bar{e}_{w\varpi}\in{\bar\A(K)}_{M}/\I_t{\bar\A(K)}_{M}=
{\A(K)}/\I_t{\A(K)}$. Hence the result.
\end{proof}
\subsection{Some results for Weyl groups}
We want to work with standard
parabolics only, and $w(P)\subset R_{Q,+}$ does not need to be
standard. We resolve this by combining terms according to left
$W_Q$ cosets. We use the following results (see \cite{C}, Section 2.7).
\begin{prop} Let $P,Q\in\P$.
The set $D^{Q,P}:=(W^Q)^{-1}\cap W^P$ intersects every double
coset $W_QwW_P$ in precisely one element $d=d(w)$, which is the
unique element of shortest length of the double coset.
\end{prop}
\begin{prop}(Kilmoyer)
Let $d\in D^{Q,P}$. Then $W_Q\cap W_{d(P)}$ is the standard
parabolic subgroup of $W_0$ corresponding to the subset
$L=Q\cap d(P)$.
\end{prop}
Let $t\in r_PT^P$ be $W_P$-generic as before, where
$W_Pr_P\subset T_P$ is the central character of a discrete
series representation $\d$ of $\H_P$. Let $\varpi=W_Pt$ be
the equivalence class of $t$.
\begin{cor}\label{cor:d}
\begin{enumerate}
\item[(i)] Let $w\in W^P$ be such that $w\varpi$ is $Q$-tempered.
Then $w(P)\subset R_{Q,+}$. We can write $w=ud$ with
$d=d(w)\in D^{Q,P}$ and $u\in W_Q$.
Then $d(P)\subset Q$, and $u\in W_Q^{d(P)}$. Conversely, if $d\in
D^{Q,P}$ is such that $d(P)\subset Q$, then for all $u\in
W_Q^{d(P)}$ we have $|ud(\varpi)|\leq_Q 1$ (in other words,
is $Q$-tempered).
\item[(ii)] The classes
$\varpi_{u,d}=:ud(\varpi)$ with $d\in D^{Q,P}$ such that
$d(P)\subset Q$ and $u\in W_Q^{d(P)}$, are mutually
disjoint.
\end{enumerate}
\end{cor}
\begin{proof}
(i) According to a result of Howlett (cf. \cite{C}, Proposition
2.7.5), we can uniquely decompose $w$ as a product of the form
$w=udv$ with $d=d(w)\in D^{Q,P}$, $u\in W_Q\cap W^L$ (with
$L=Q\cap d(P)$), and $v\in W_P$.
In fact $v=e$, since otherwise there would exist a $\a\in
R_{P,+}$ with $v(\a)=-\a_p\in P$. But then $ud(\a_p)<0$, which
implies according to \cite{C}, Lemma 2.7.1 that $d(\a_p)=\a_q\in
L$. Hence $u(\a_q)<0$, which contradicts the assumption $u\in
W_Q\cap W^L$. Thus we have $d(P)=u^{-1}w(P)\subset R_{Q,+}$, whence
$W_{d(P)}\subset W_Q$. By Kilmoyer's result it now follows
that $W_{d(P)}=W_{Q\cap d(P)}$. Hence $d(P)\subset Q$ and $L=d(P)$.
The converse is clear.

(ii) Suppose that
$\varpi_{u,d}\cap \varpi_{u^\prime,d^\prime}\not=\emptyset$.
The Weyl
group $W_0$ permutes equivalence classes, thus this implies
that $(ud)^{-1}u^\prime d^\prime(t)\in \varpi$. Since $t$ is
generic, there exists a $v\in W_P$ such that
$u^\prime d^\prime=udv$. By Howlett's result \cite{C},
Proposition 2.7.5 this implies that $v=1$, $u=u^\prime$
and $d=d^\prime$.
\end{proof}
\begin{cor}\label{cor:def}
For all $d\in D^{Q,P}$ with $d(P)\subset Q$
we write
\[e_{W_Qd\varpi}=\sum_{u\in W_Q^{d(P)}} e_{ud\varpi}.\]
This is a collection of orthogonal idempotents
of $\bar\H_t$. The constant term of $f=f_{a,b}(\xi)$
equals
\[
f^Q(h)=\sum_{d\in D^{Q,P}:d(P)\subset Q} f^d(h),
\]
where we define $f^d(h):=f(e_{W_Qd\varpi}h)$. This is the
contribution to the constant term $f^Q$ of $f$ whose exponents
have the same restriction to $X\cap Q^\perp$ as $d(t)$.
\end{cor}
\subsection{The singularities of $f^d$}
{\it In this section we take the formulae of Corollary \ref{cor:def}
as a definition
of $f^Q$ and $f^d$, even when $t^P\in T^P$ is not in $T^P_u$.}

We will now bound the possible singularities of the individual
contributions $f^d$ to $f^Q$, viewed as functions of $t^P\in T^P$.
We have seen in Corollary \ref{cor:rat} that
$f^d$ extends to a rational function of
$\xi\in\Xi$. To stress this dependence we sometimes
write $f^d(\xi,h)$. We write $\xi=(P,\d,t^P)$ and
put $t=t(\xi)=r_Pt^P$, where $r_P\in T_P$
is such that $W_Pr_P$ is equal to the central
character of $\d$.
\begin{lem}\label{lem:basic}
Let $P,Q\in\P$ and let $d\in D^{Q,P}$ be such that
$d(P)\subset Q$.
Let $h,h^\prime\in\bar{\H}_t$. Then
\begin{equation}
f^d_{a,b}(\xi;hh^\prime)=
f^d_{a,\pi(\xi;h^\prime)(b)}(\xi;h).
\end{equation}
\end{lem}
\begin{proof}
This follows immediately from Corollary \ref{cor:basic}.
\end{proof}
\begin{lem}\label{lem:basic2}
As in Lemma \ref{lem:basic}.
Let $g\in\W_P$ and put $P^\prime=g(P)$.
According to Corollary \ref{cor:d} we can
write $dg^{-1}=u^\prime d^\prime$ with
$d^\prime\in D^{Q,P^\prime}$ and
$u^\prime\in W_Q^{P^\prime}$. We put $t^\prime=g(t)$ and
$\varpi^\prime=W_{P^\prime}t^\prime=g(\varpi)$,
so that $e_{W_Qd\varpi}=e_{W_Qd^\prime\varpi^\prime}$.
With these notations we have
the following equality of rational
functions of $\xi$:
\begin{equation}\label{eq:trans}
f_{a,b}^d(\xi;h)=
f_{\pi(g,\bar{\xi}^{-1})(a),\pi(g,\xi)(b)}^{d^\prime}
(g(\xi);h),
\end{equation}
where $\bar{\xi}^{-1}:=(P,\d,\bar{t^P}^{-1})$.
\end{lem}
\begin{proof}
This equation follows from the special case $\xi\in\Xi_{P,\d,u}$
because the left hand side and the right hand side are obviously
rational functions of $\xi$. In this special case the equation
simply expresses the unitarity of the intertwiners (cf. Theorem
\ref{thm:fusm}).
\end{proof}
\begin{lem}
Let $P,Q\in \P$. Then $\H$ has the following
direct sum decomposition in
left $\H^Q$-right $\H(W_P)$-submodules:
\begin{equation}
\H=\bigoplus_{d\in D^{Q,P}}\H_{Q,P}(d),
\end{equation}
where $\H_{Q,P}(d):=\H^QN_d\H(W_P)$.
\end{lem}
\begin{proof}
Using the Bernstein presentation of $\H^Q$ and the
definition of the multiplication in $\H(W_0)$ we
easily see that
\begin{equation}\label{eq:hqp}
\H_{Q,P}(d)=\bigoplus_{w\in W_QdW_P}\A N_w.
\end{equation}
The result thus follows from the
Bernstein presentation of $\H$.
\end{proof}
After these preparations we will now concentrate on
an important special case.
\begin{dfn}
Let $\pi_{Q,P}^d:\H\to\H_{Q,P}(d)$ denote the projection
according to the above direct sum decomposition.
Given $Q\in\P$, denote by
$w^Q=w_0w_Q$
\index{wa@$w_0$, longest element of $W_0$}
\index{wb@$w_P$, longest element of $W_P$}
\index{wc@$w^P$, longest element of $W^P$}
the longest element
of $W^Q$, and by $Q^\prime=w^Q(Q)=-w_0(Q)\in\P$. Then
$w^{Q^\prime}=(w^Q)^{-1}\in D^{Q,Q^\prime}$, and
\begin{equation}\label{eq:dfn}
\H_{Q,Q^\prime}(w^{Q^\prime})=\H^QN_{w^{Q^\prime}}
=\A N_{w^{Q^\prime}}\H(W_{Q^{\prime}}).
\end{equation}
Let $p_{Q}:\H\to\H^Q$ be the left
$\H^Q$-module map defined by
\begin{equation}
p_Q(h):=\pi_{Q,Q^\prime}^{w^{Q^\prime}}(h)N_{w^{Q^{\prime}}}^{-1}
\end{equation}
(Observe that this map indeed has values in $\H^Q$ by
(\ref{eq:hqp})).
\end{dfn}
In (\ref{eq:dfn}) we have used that $N_{w^{Q^\prime}}N_{w^\prime}=
N_wN_{w^{Q^\prime}}$ if $w\in W_Q$.
\begin{thm}\label{thm:spec}
Let $P,Q\in \P$ be such that $P\subset Q$. We put
$P^\prime:= w^Q(P)\subset Q^\prime\in\P$ and
$\xi^\prime=w^Q(\xi):=(P^\prime,\d^\prime,t^{P^\prime})$.

Let $a^\prime\in i(V_{\d^\prime})=\H(W^{P^\prime})\otimes V_{\d^\prime}$,
$b^\prime\in\H(W_{Q^\prime}^{P^\prime})\otimes V_{\d^\prime}\subset
i(V_{\d^\prime})$ and let $h\in \H$.
We introduce the unitary isomorphism
\begin{equation}
\s:=\psi_{w^Q}\otimes\tilde{\delta}_{w^Q}:
\H(W_Q^P)\otimes V_\d\to
\H(W_{Q^\prime}^{P^\prime})\otimes V_{\d^\prime},
\end{equation}
and the orthogonal projection
\begin{equation}
\rho:i(V_\d)\to \H(W^P_Q)\otimes V_\d.
\end{equation}
With these notations, put
\begin{align}
a:&=\rho(\pi(w^Q,\bar{\xi}^{-1})^{-1}(a^\prime))
\in\H(W^P_Q)\otimes V_\d\\
\nonumber
b:&=\s^{-1}(b^\prime)\in\H(W^P_Q)\otimes V_\d
\end{align}
We then have, with
$c^Q(\xi):=\prod_{\a\in R_{0,+}\backslash R_{Q,+}}c_\a(t)$,
\begin{equation}
f^{w^{Q^\prime}}_{a^\prime,b^\prime}(\xi^\prime,h)
=q(w^{Q})^{1/2}c^Q(\xi)f_{Q,a,b}(\xi,p_Q(h)).
\end{equation}
Here $f_{Q,a,b}(\xi,h)=f_{a,b}(\xi,h)$
(with $h\in\H^Q$, $a,b\in\H(W^P_Q)\otimes V_\d$)
is the matrix coefficient (associated to the pair
$a,b$) of the representation
\begin{equation}
\pi^Q(\xi):=\operatorname{Ind}_{\H^P}^{\H^Q}\d_{t^P}
\end{equation}
of $\H^Q$ (which is tempered and unitary if $\xi\in \Xi_u$).
\end{thm}
\begin{proof}
Choose $r_P\in T_P$ such that $W_Pr_P$ is the central character
of $\d$, and write $t^\prime=w^Q(t)$ with $t=r_Pt^P$.
Since we are dealing with rational functions of $\xi$
it is sufficient to assume that $\xi$ is regular, i.e.
that $t$ is $R_P$-regular. We then extend
$\pi(\xi^\prime)$ to
the completion $\bar{\H}_t$ (recall \ref{subsub:ext})
and study $\pi(\xi^\prime)$
in the light of the isomorphisms (\ref{eq:dec}) and (\ref{eq:iso}).

We combine, in the decomposition (\ref{eq:dec})
applied to the parabolic $P^\prime=w^Q(P)$ and parameter $t^\prime$,
the idempotents according to left cosets of $W_{Q^\prime}$.
In other words, we partition $W_0t$ into the sets
$w(\Omega)$ with $w\in W^{Q^\prime}$ and
$\Omega=W_{Q^\prime}t^\prime=W_{Q^\prime}^{P^\prime}\varpi^\prime$
(with $\varpi^\prime=w^Q(\varpi)=W_{P^\prime}t^\prime$). These
sets are evidently unions of the original equivalence classes
in formula (\ref{eq:dec}) (with respect to $P^\prime$ and
$t^\prime$),
the left $W_{P^\prime}$-cosets acting on $t^\prime$.
We denote the corresponding idempotents by
(for all $w\in W^{Q^\prime}$)
\[
e^\sharp_w:=\sum_{x\in W^{w^Q(P)}_{Q^\prime}}e_{wx\varpi^\prime}.
\]

Note that $t^\prime$ is $P^\prime$-generic, and thus certainly
$Q^\prime$-generic. The structure formula (\ref{eq:dec})
holds therefore, also in terms of the idempotents $e^\sharp_w$,
where we replace in (\ref{eq:dec}) the parabolic $P^\prime$
by $Q^\prime$.

Remark that $e^\sharp_{w^{Q^\prime}\Omega} N_w e^\sharp_\Omega=0$
for any $w\in W_0$ with length of less than $|R_{0,+}\backslash
R_Q|$ (=the length of $w^{Q^\prime}$). Note by the way that
$e^\sharp_{w^{Q^\prime}\Omega}=e^\sharp_{W_Qt}$.
Thus for all $d\in D^{Q,Q^\prime}$, $d\not= w^{Q^\prime}$
we see that
$e^\sharp_{W_Qt} \H^{Q,P}(d) e^\sharp_{w^QW_Qt}=0$.

Hence for all $h\in\H$, $a^\prime\in i(V_{\d^\prime})$ and
$b^\prime\in\H(W^{P^\prime}_{Q^\prime})$ we have
\begin{equation}
f^{w^{Q^\prime}}_{a^\prime,b^\prime}
(\xi^\prime,h) =f_{a^\prime,b^\prime}
(\xi^\prime,p_Q(h) e^\sharp_{W_Qt} N_{w^{Q^\prime}}
e^\sharp_{w^QW_Qt}).
\end{equation}
Since $f_{a^\prime,b^\prime}^{w^{Q^\prime}}(\xi^\prime,\H\I_t)=0$
we can use the analog of formula (4.57) of \cite{O}
(we use here that the c-function $c^Q(t)$ is
$W_Q$-invariant, together with the argument in the proof of the
Proposition \ref{cor:idmp}. This makes that we can evaluate the
c-factors at $t^\prime$):
\begin{equation}
f^{w^{Q^\prime}}_{a^\prime,b^\prime}
(\xi^\prime,h)=
q(w^Q)^{1/2}c^Q(t) f_{a^\prime,b^\prime} (\xi^\prime,p_Q(h)
\iota^0_{w^{Q^\prime}}).
\end{equation}
We use Lemma \ref{lem:basic} and then rewrite the result
using (\ref{eq:dec}) and Definition \ref{dfn:A}. Assume that
$b^\prime=x^\prime\otimes v^\prime$ and
$b=\s^{-1}(b^\prime)=x\otimes v$. Then
\begin{align*}
\iota^0_{w^{Q^\prime}}(b^\prime)&=
\iota^0_{w^{Q^\prime}}
(x^\prime e^\sharp_{W_{Q^\prime}t^\prime}\otimes v^\prime)\\
&=e^\sharp_{W_Qt}(\psi_{w^{Q^\prime}}(x^\prime)
\iota^0_{w^{Q^\prime}}\otimes v^\prime)\\
&=e^\sharp_{W_Qt}\pi(w^Q,\xi)(x\otimes v)\\
&=\pi(w^Q,\xi)(b).
\end{align*}
Thus we obtain
\begin{align}
f^{w^{Q^\prime}}_{a^\prime,b^\prime}(\xi^\prime,h)&=
q(w^Q)^{1/2}c^Q(\xi)f_{a^\prime,\pi(w^Q,\xi)(b)}(\xi^\prime,p_Q(h))\\
\nonumber&=
q(w^Q)^{1/2}c^Q(\xi)f_{(\pi(w^Q,\bar{\xi}^{-1})^{-1}
(a^\prime)),b}(\xi,p_Q(h))\\
\nonumber&=
q(w^Q)^{1/2}c^Q(\xi)f_{Q,a,b}(\xi,p_Q(h)).
\end{align}
In the second step we use the unitarity of the
intertwining operators $\pi(w^Q,\xi)$
to rewrite the matrix coefficient as a coefficient
of the induced representation $\pi(\xi)$ (extended
holomorphically as in Lemma \ref{lem:basic2}; in fact it
is a simple special case of this Lemma).
Since $b\in \H(W^P_Q)\otimes V_\d$ and $p_Q(h)\in\H^Q$,
we can project the vector $\pi(w^Q,\bar{\xi}^{-1})^{-1}(a^\prime)$
onto $\H(W^P_Q)\otimes V_\d$, and consider the result as
a matrix coefficient of $\pi^Q(\xi)$.
\end{proof}
\begin{thm}\label{thm:cor}
Fix $P\in \P$ and $\d\in\Delta_{P,W_{P}r_P}$.
We denote by $\Xi_{P,\d}\subset\Xi$ the collection of standard
induction data of the form $(P,\d,t^P)$ with $t^P\in T^P$, and
by $\Xi_{P,\d,u}\subset\Xi_{P,\d}$ the subset of such triples
with $t^P\in T^P_u$.
Then for all $d\in D^{Q,P}$ such that $d(P)\subset Q$
and for all $a,b\in i(V_\d)$, the
rational function
\begin{equation}
\xi\to c(\xi)^{-1}f^d_{a,b}(\xi,h)
\end{equation}
is regular in a neighborhood of $\Xi_{P,\d,u}$.
\end{thm}
\begin{proof}
We apply Lemma \ref{lem:basic2} with $g=w^Qd\in W(P,P^\prime)$
where $P^\prime=w^Q(d(P))\subset Q^\prime$. Notice that
$d^\prime=w^{Q^\prime}$. Put $\xi^\prime=g(\xi)$ and
\begin{align}
a^\prime&=\pi(g,\bar{\xi}^{-1})(a)\\
\nonumber
\tilde{b}^\prime&=\pi(g,\xi)(b).
\end{align}
We obtain
\begin{equation}
c(\xi)^{-1}f_{a,b}^d(\xi,h)=
c(\xi)^{-1}
f_{a^\prime,\tilde{b}^\prime}^{w^{Q^\prime}}(\xi^\prime,h).
\end{equation}
Now we can {\it uniquely} decompose $\tilde{b}^\prime$ in the
following way
\begin{equation}
\tilde{b}^\prime=\pi(\xi^\prime,\tilde{h})(b^\prime)
\end{equation}
with $\tilde{h}\in\H(W^{Q^\prime})$ and
$b^\prime\in\H(W_{Q^\prime}^{P^\prime})\otimes V_{\d^\prime}\subset
i(V_{\d^\prime})$.
With the help of Lemma \ref{lem:basic} we get
\begin{equation}
c(\xi)^{-1}f_{a,b}^d(\xi,h)=
c(\xi)^{-1}
f_{a^\prime,b^\prime}^{w^{Q^\prime}}(\xi^\prime,h\tilde{h}).
\end{equation}
We can now apply Theorem \ref{thm:spec} with
respect to $P_0:=d(P)\subset Q$. We put $\xi_0:=d(\xi)$
and
\begin{align}
a_0:&=\rho(\pi(w^Q,\bar{\xi_0}^{-1})^{-1}(a^\prime))
\in\H(W^{P_0}_Q)\otimes V_\d\\
\nonumber
b_0:&=\s^{-1}(b^\prime)\in\H(W^{P_0}_Q)\otimes V_\d
\end{align}
to obtain:
\begin{equation}
c(\xi)^{-1}f_{a,b}^d(\xi,h)=
q(w^Q)^{1/2}(c(\xi)^{-1}c(\xi_0))c_Q(\xi_0)^{-1}
f_{Q,a_0,b_0}(\xi_0,p_Q(h\tilde{h})),
\end{equation}
where in general for $Q\supset P$ we denote
\begin{equation}\label{eq:rel}
c_Q(\xi):=\prod_{\a\in R_{Q,+}\backslash R_{P,+}}c_\a(t).
\end{equation}
The regularity of the normalization factor $c_Q(\xi_0)^{-1}$
as a function of $\xi_0$ (and thus as a function of
$\xi=d^{-1}(\xi_0)$) follows from \cite{O}, Theorem 3.25,
when we consider the tempered residual coset $r_PT^P_u\subset T$
for the Hecke algebra $\H^Q$ (instead of $\H$ itself). It is a simple
special case of Proposition \ref{lem:csmo}(v). Similarly,
the regularity of $c(\xi)^{-1}c(\xi_0)$ is asserted by
Proposition \ref{lem:csmo}(iv).
By the regularity of the various intertwining operators we have
used (cf. Theorem \ref{thm:fusm}) it is clear that also $a_0, b_0$
are also rational and regular on $\Xi_{P,\d,u}$. We have finished
the proof.
\end{proof}
From the proof of Theorem \ref{thm:cor} we obtain
\begin{cor}\label{cor:fina}
Thus, the final conclusion of all these considerations is that
for all $h\in\H$, $P\in\P$, $\d\in\Delta_P$,
and $a,b\in i(V_\d)$ fixed, the function
\begin{equation}
\Xi_{P,\d,u}\times\H^Q\ni(\xi,h^Q)\to c(\xi)^{-1}f_{a,b}^Q(\xi,h^Qh)
\end{equation}
is a linear combination
with coefficient functions which are regular rational functions
on $\Xi_{P,\d,u}$, of normalized matrix coefficients
\begin{equation}
c_Q(d(\xi))^{-1}f_{Q,a^\prime,b^\prime}(d(\xi),h^Q)
\end{equation}
of induced representations of $\H^Q$ of the form
$\pi^Q(d(\xi))$ (where $d$ ranges over the Weyl group elements
$d\in D^{Q,P}$ such that $d(P)\subset Q$).
\end{cor}
\section{Proof of the main theorem}
\subsection{Uniform estimates for the  coefficients of $\pi(\xi)$}
Recall that $X^+$ is the cone
$\{x\in X\mid \langle x,\a^\vee\rangle\geq 0
\ \mathrm{for\ all\ }\a\in R_{1,+}\}$.
We put $Z_X=X^+\cap X^-$. This is a sublattice of elements in
$X$ with length $0$. Recall that $Q$ denotes the root lattice. The
sublattice
$Q\oplus Z_X\subset X$ has finite index in $X$. If $x=x_Q+x_Z\in
Q^++Z_X$ then
\begin{equation}\label{eq:n}
\mathcal{N}(x)=x_Q(2\rho^\vee)+\Vert x_Z\Vert.
\end{equation}
Let us show that $Q^{+}$ is finitely generated over $\Z_{+}$.
For each fundamental weight $\delta_{i}$, let
$q_i=m_i\delta_i$ be the least multiple of $\d_i$ such
that $q_i\in Q$ (thus
$m_i\in\N=\{1,2,3,\cdots\}$
is a divisor of the index
$[P:Q]$). These multiples generate over
$\Z_+$ a cone
$C^+\subset Q^+$. Let $F$ be the fundamental
domain $F=\{\sum_i t_iq_i |t_i\in [0,1)\}$ of $C$, and let
$F_Q=F\cap Q\subset Q^+$ (a finite set).
Clearly $F_Q$ and the
$\{q_i\}$ generate $Q^+$ over $\Z_+$.
Let $x_1,\cdots,x_m,x_{m+1},\cdots,x_N\in X^+$
such that $x_1,\cdots,x_m$ is a set of $\mathbb{Z}_+$-generators
of $Q^+$ and that $x_{m+1},\cdots,x_N\in Z_X$
is a $\mathbb{Z}$-basis of $Z_X$. By (\ref{eq:n}) we see that
there exists a constant $K>0$ such that for all
$x\in Q^++X_Z$ and all decompositions
$x=\sum l_ix_i$ with $l_i\geq 0$ if $i\leq m$, we have
\begin{equation}\label{eq:est}
\sum|l_i|\leq K\mathcal{N}(x)
\end{equation}
(just observe that $x_i(2\rho^\vee)\geq 1$ if $i\in 1\dots,m$).
We fix such a $K>0$.

We define a function $\nu$ on $T_{rs}$ by
\begin{align}
\nu(t)=
\operatorname{max}
(\nonumber\{|x_i(t)|\mid &i=m+1,\dots,N\}\\
&\cup\{\vert x_i(wt)\vert \mid i=1,\dots,m;w\in W_0\}).
\end{align}
\index{0n@$\nu$, (multiplicative) distance function on $T_{rs}$}
The positive real cone spanned by the elements
$wx_i$ ($w\in W_0$,
$i=1,\dots,m$) and $\pm x_i$ ($i=m+1,\dots,N$) is
the full dual of $\operatorname{Lie}(T_{rs})$.
Therefore the function $\log(\nu)\circ\exp$ is
a norm on $\operatorname{Lie}(T_{rs})$.
\begin{thm}\label{thm:rough}
Let $R>1$, $P\in\P$, and $\d\in\Delta_P$ be given.
Choose a set $x_i\in X^+$ as above and let
$K$ and $\nu$ be as above. We use the notation
$\nu(|\xi|):=\nu(|t^P|)$ for
$\xi=(P,\d,t^P)\in\Xi_{P,\d}$.
Define a compact neighborhood $D^P(R)\subset \Xi_{P,\d}$
of $\Xi_{P,\d,u}\subset \Xi_{P,\d}$ by
$D^P(R)=\{\xi\in\Xi_{P,\d}\mid \nu(|\xi|)\leq R\}$.

There exists a $d\in\N$, and
there exists a constant $c>0$
(depending on $R$ only) such that
for all $w\in W$, for all $a,b\in i(V_\d)$,
and for all $\xi\in D^P(R)$,
the matrix coefficient $f_{a,b}(\xi,N_w)$
satisfies
\begin{equation}\label{eq:rough}
|f_{a,b}(\xi,N_w)|
\leq c\Vert a\Vert\Vert b\Vert
(1+\mathcal{N}(w))^d\nu(|\xi|)^{K\mathcal{N}(w)}
\end{equation}
\end{thm}
\begin{proof}
Using
\cite{O}, equation (2.27) (also see the proof of Lemma
\ref{lem:w=uxv})
we see that it is equivalent to
show that $f_{a,b}(\xi,N_u\theta_x N_v)$ can be estimated
by the right hand side of (\ref{eq:rough}) with $w=x$, for
all $u,v\in W_0$ and $x\in X^+$.
By applying right (resp. left) translations of the matrix coefficient
$f_{a,b}(\xi)$ by $N_v$ (resp. $N_u$) and by a set of representatives
of the finite quotient $X/(Q+ Z_X)$ we see that we may further
reduce to proving the estimates (\ref{eq:rough}) for
$w=x\in Q^++ Z_X$.

Recall (cf. \cite{O}, Proposition 4.20 and its proof)
that the the eigenvalues of  the matrix of
$\pi(\xi,\theta_x)$  are of the form $x(w_i(r_{j}t^P))$.
Here the $r_{j}\in T_P$ are the generalized $X_P$-eigenvalues
of the discrete series representation $\d$. By Casselman's
criterion we know therefore that for all $x\in X^+$,
$x(w_i(r_{j}))\leq 1$. This implies that for all  $i\in 1,\dots,m$, and
for all $\xi\in \Xi_{P,\d}$, the  eigenvalues of
$\pi(\xi,\theta_{x_{i}})$ are bounded by
$\nu(|\xi|)$.
Then Lemma \ref{lem:ban} allows to estimate the norm of
$\pi(\xi,\theta_{l_{i}x_{i}})$, by dividing
$\pi(\xi, \theta_{x_{i}})$ by
$\nu (\vert \xi \vert )$.

Taking into account the fact that $D^{R}(P)$ is compact,
one sees that the norm of
$\pi(\xi,\theta_{x_{i}})$ is bounded if $\xi$ is in   $D^{R}(P)$.
By a simple product formula, one estimates the norm of
$\pi(\xi,\theta_{x})$.
These estimates together with equation (\ref{eq:est})
imply the desired result.
\end{proof}
\begin{cor}\label{cor:derest}
For all constant
coefficient differential operators $D$ on $\Xi_{P,\d}$
there exist constants $d\in \N$ and $c>0$
such that for all $\xi\in\Xi_{P,\d,u}$, for all $a,b\in i(V_\d)$,
and for all $w\in W$
\begin{equation}
|Df_{a,b}(\xi,N_w)|\leq c\Vert a\Vert\Vert b\Vert
(1+\mathcal{N}(w))^d.
\end{equation}
\end{cor}
\begin{proof}
This is a standard application of the Cauchy integral formula,
starting with equation (\ref{eq:rough}).
Choose a basis $x_1,\dots,x_p$ of the character lattice
$X^P$ of $T^P$, and let $y_1,\dots,y_p$ be the dual basis.
Let $C_\epsilon:=
\{v\in\operatorname{Lie}(T^P)\mid\forall i:\ |x_i(v)|=\epsilon\}$.
We may assume that $D$ is of the form
$D=D^\a:=y_1^{\a_1}\dots y_p^{\a_p}$.
By the holomorphicity of $f_{a,b}$ we have,
for a suitable constant $C_\a>0$ and any
choice of a sequence of radii $\epsilon(w)$:
\begin{equation}\label{eq:cauchy}
D^\a f_{a,b}(\xi,N_w)=C_\a\int_{v\in C_{\epsilon(w)}}
\frac{f_{a,b}(\exp(v)\cdot\xi,N_w)}{\prod_i x_i(v)^{{\a}_i+1}}
dx_1\wedge\dots\wedge dx_p.
\end{equation}
Now use the estimates of Theorem \ref{thm:rough} with the
sequence $\epsilon(w)$ chosen such that
$r(w):=\operatorname{max}\{\nu(|\exp(v)|)\mid v\in C_{\epsilon(w)}\}$
is equal to
\begin{equation}\label{eq:a1}
1+1/(1+\mathcal{N}(w)).
\end{equation}
But $\log(\nu)\circ\exp$ is a norm on $\operatorname{Lie}(T^P)$,
as well as
$\operatorname{Sup}_{i}\vert x_i(v)\vert$.
They are equivalent. Moreover log$(1+x)\geq k' x$ for $x\in [0,1]$,
for
some $k'>0$. Together with \ref{eq:a1}
this implies that there exists a
constant
$k>0$ such that
$\epsilon(w)\geq k/(1+\mathcal{N}(w))$. So equation
(\ref{eq:cauchy})
yields the estimate (for some constant $c^\prime>0$)
\begin{equation}
|D^\a f_{a,b}(\xi,N_w)|\leq c^\prime\Vert a\Vert\Vert b\Vert
(1+\mathcal{N}(w))^{d+|\a|}
(1+1/(1+\mathcal{N}(w)))^{K\mathcal{N}(w)}.
\end{equation}
This easily leads to the desired result.
\end{proof}
\begin{cor}\label{cor:imF}
We have
$\F(\S)\subset \cc(\Xi_u,\operatorname{End}(\V_{\Xi}))^\W$.
The restriction $\F_\S$ of $\F$ to $\S$ defines a continuous
map $\F_\S:\S\to\cc(\Xi_u,\operatorname{End}(\V_{\Xi}))^\W$.
\end{cor}
\begin{proof}
The equivariance of the sections in the image is clear.
Recall that
$\F(N_w)\in
\operatorname{Pol}(\Xi_u,\operatorname{End}(\V_{\Xi}))$
is defined by
$\F(N_w)(\xi)=\pi(\xi,N_w)$.

Hence by the estimates of Corollary \ref{cor:derest} we see that
for any continuous seminorm
$p$ on $\cc(\Xi_u,\operatorname{End}(\V_{\Xi}))$ there exist
constants $C>0$ and
$d\in\Z_+$
such that
$p(\F(N_w))\leq C(1+\mathcal{N}(w))^d$.

Let $b\in\Z_+$ be such that
$0<C_b:=\sum_{w\in W}(1+\mathcal{N}(w))^{-b}<\infty$, and let
$q=q_p$ denote the continuous seminorm on $\S$ defined by
$q(x):=CC_b\operatorname{sup}_w\vert (x,N_w)
\vert(1+\mathcal{N}(w))^{d+b}$.
Then $p(\F_\S(x))\leq q(x)$ for all $x\in\H$, implying that
$\F_\S$ is a continuous map as claimed.
\end{proof}
\subsection{Uniform estimates of the difference of a coefficient
and its constant term}
We now introduce the important notion of a normalized smooth
family of coefficients:
\begin{dfn}
Let $P\in {\mathcal P}$ and let $\delta\in \Delta_P$ be an
irreducible discrete series of ${\mathcal H}_{P}$ with central character
$W_Pr_P\in W_P\backslash T_P$.
We put $\xi=(P,\d,t^P)\in \Xi_{P,\d,u}$.
A smooth family of coefficients of
$\pi(\xi)$, $\xi\in\Xi_{P,\d,u}$ is a family of linear functionals
on $\H$ of the form:
\begin{equation}\label{eq:fam}
\H\ni h\mapsto
\operatorname{Tr}(\sigma(\xi)\pi(\xi)(h)),
\end{equation}
where $\sigma$ is a section of
$\cc(\Xi_{P,\d,u},\operatorname{End}(\V_\Xi))$.

A smooth section $\sigma\in\cc(\Xi_{P,\d,u},\operatorname{End}(\V_\Xi))$
is called normalized smooth when it is divisible (in the
$\cc(\Xi_{P,\d,u})$-module
$\cc(\Xi_{P,\d,u},\operatorname{End}(\V_\Xi))$)\index{C@$\cc(\Xi_u)$,
space of $\cc$-functions on
$\Xi_u$}
by the smooth function
$\{\xi\to c^{-1}(\xi)\}\in\cc(\Xi_{P,\d,u})$
(cf. Proposition \ref{lem:csmo}).

A normalized smooth family of coefficients of
$\pi(\xi)$, $\xi\in\Xi_{P,\d,u}$ is a smooth family of
coefficients (\ref{eq:fam}) for which $\sigma$ is
normalized smooth.
\end{dfn}
\begin{rem}
We frequently use $t^P$ rather than $\xi=(P,\d,t^P)$ as the
parameter of a family of coefficients.
\end{rem}
\begin{cor}\label{cor:trans}
It follows directly from the definitions that smooth
(resp. normalized smooth) families of coefficients of
$\pi(\xi)$, $\xi\in\Xi_{P,\d,u}$, are stable under left and right
translations by elements of $h\in \H$.
\end{cor}
Moreover, Corollary \ref{cor:fina} implies that:
\begin{cor}
The restriction to ${\mathcal H}^{Q}$  of the constant term of a
normalized smooth family of
$\pi(\xi)=\pi(P,\delta,t^{P})$, $t^{P}\in T^{P}_{u}$, along $Q\in
{\mathcal P}$ is a finite sum of terms, each of these being a
normalized smooth
family of coefficients of $\pi^{Q}(d(\xi))$, where $d$ is some
Weyl group
element with $d(P)\subset Q$.
\end{cor}
\begin{lem}\label{lem:est1}
Assume $Z_{X}=\{0\}$.
Let $\Xi_{P,\d,t^P}\ni\xi\to\Phi_\xi$
be a smooth family of coefficients.

Let $\a\in F_{0}$, and put $Q=F_{0}\setminus\{\alpha\}$.

Let $\Vert\cdot\Vert$ denote a norm which comes from a
$W_0$-invariant euclidean structure on $X\otimes_{\Z}\R$.

Let $a>0$ and let $X^+_{a}$ denote the cone (over $\Z_+$)
$X^+ _{a}=\{x\in X^{+}\mid \langle x,\a^\vee\rangle
>a\Vert  x\Vert\}$.
Then there exists $C, b>0$ such that
\begin{equation}\label{eq:est1}
\vert (\Phi_\xi- \Phi_\xi^{Q})(N_u\theta_{x}N_v)
\vert\leq C e^{-b \Vert x\Vert },
\end{equation}
for all $x\in X^{+}_{a}$, $\xi\in\Xi_{P,\d,u}$, $u\in W_Q$, and
$v\in W_{0}$.
\end{lem}
\begin{proof}
Recall that the lattice $X$ contains the root lattice $Q(R_0)$,
and hence an
integral multiple of the weight lattice $P(R_0)$, $kP(R_0)$ say.
We put $X^\prime=kP(R_0)\subset X$ and we identify $X^\prime$
with $\Z^{l}$ via a basis of $X^\prime$ consisting of the elements
$(k\d _{\b}), \>\>\b\in F_{0}$ (where the $\delta_{\b}$ are the
fundamental weights), ordered in such a way that $e_{1}=k\d_{\a}$.
The temperedness of $\pi(P,\delta,t^{P})$, $t^{P}\in
T^{P}_{u}$, and the fact that its central character is given by
$t=r_{P}t^{P}$ imply that the possible eigenvalues of
$\pi(P,\delta,t^{P})(\theta_{k\delta_{\b}})$  are among the
$wt(k\d_\b)$ with $w\in W_{0}$ such that
$\vert wt(k\delta_{\b})\vert \leq 1$.
Moreover the modulus of $wt(k\delta_{\b})$, hence of
$wt(e_{1})$, does not depend on $t^{P}\in T_{u}^{P}$.

Hence if  $u=v=e$ and $x\in X'=kP(R_0)$, (\ref{eq:est1}) follows,
in view of the definition of the constant term,
from Lemma \ref{lem:specproj}.

Let us now derive the general case of (\ref{eq:est1}) from this
special case. Since $Z_{X}=\{0\}$, $X'$ is of finite index in $X$.
One can assume that
$a$ is small enough, in such a way that $X^{+}_{a}$ is nonempty,
otherwise there is nothing to prove.
Let $x\in X^+_a$ and let $y$ be the orthogonal projection
of $x$ on the line $\R\delta_{\alpha}$. By definition of $X_Q$
we have $x-y\in X_Q$, and since
$\langle x-y,\b^\vee\rangle=\langle x,\b^\vee\rangle$ for all
$\b\in Q$ we find that in fact $x-y\in X_Q^+$. Thus $x-y$ is
a nonnegative
linear combination of the fundamental weights $\d^Q_\b$ of $R_Q$.
It is a basic fact that the fundamental weights
of a root system (with given basis of simple roots)
have nonnegative rational coefficients in the basis of
its simple roots (indeed, this statement reduces to the case of an
irreducible root system, in which case the indecomposability
of the Cartan matrix implies that these coefficients
are in fact strictly positive). Hence we have
$x-y\in X_Q^+\subset \Q_+Q$. Since $\langle\b,\a^\vee\rangle\leq 0$
for all $\b\in Q$, we see that
$\langle y,\a^\vee\rangle\geq\langle x,\a^\vee\rangle>
a\Vert x\Vert\geq a\Vert  y\Vert$. Hence
$k\delta_{\alpha}\in  X^{+}_{a}$.

Let
$(x_{1}, \dots,x_{r})$ be a set of representatives in $X$ of
$X/X'$. Let us show that one can choose the $x_{i}$ in
$-X^{+}_{a}$. Our claim is a consequence of the following fact. If
$y\in X$, one has
$y+n \delta_{\alpha} \in X^{+}_{a}$ for $n$ large. In fact by the
triangular inequality one has:
$$\langle y+n \delta_{\alpha},\a^\vee\rangle -
a \Vert y+ n \delta_{\alpha}\Vert\geq
n(\langle\delta_{\alpha},\a^\vee\rangle-a\Vert\delta_{\alpha}\Vert)+\alpha(y)
-a\Vert y\Vert $$
Thus, if $x\in
X^{+}_{a}$ and $ x=x'+x_{i}$, for some $x'\in X'$ and some $i$,
one has $x'\in {X'}^{+}_{a}$. To get the estimates, one applies the
previous estimates
to the translates of the family
$\Phi_{\xi}$ by the $N_u$ (from the left),
and by $\theta_{x_{i}}N_v$ (from the right),
which are smooth families
of coefficients themselves (cf. Corollary \ref{cor:trans}),
and taking into account Corollary
\ref{cor:basic}.
\end{proof}
\subsection{Wave packets}
Recall that $\J$ was introduced as the adjoint of $\F$. Thus
if $\sigma\in L_2(\Xi_{P,\d,u},\operatorname{End}(\V_\Xi),\mu_{Pl})$ then
$\J(\sigma)\in L_2(\H)$, and is completely characterized by
the value of $(\J(\sigma),h)$ where $h\in\H$. We have, using
Theorem \ref{thm:mainO}, that
\begin{align}\label{eq:wave}
\J(\sigma)(h):&=(\J(\sigma)^*,h)=(h^*,\J(\sigma))=\\
\nonumber &=(\F(h^*),\sigma)=\mu_{\Ri,\d}\int_{\Xi_{P,\d,u}}
\operatorname{Tr}(\sigma(\xi)\pi(\xi,h))|c(\xi)|^{-2}d\xi,
\end{align}
where $\mu_{\Ri,\d}=q(w^P)^{-1}|\W_P/\K_P|^{-1}
\mu_{\Ri_P,Pl}(\{\d\})>0$.

Recall Definition \ref{dfn:calC}.
Assume that $\tilde{\sigma}=c(w^P\cdot)\sigma
\in\mathcal{C}(\Xi_u,\operatorname{End}(\V_\Xi))$
(in other words,
$\sigma\in\cc(\Xi_{P,\d,u},\operatorname{End}(\V_\Xi))$).

Denote by $\Phi^{\sigma}$ the smooth family of
coefficients $\Phi_\xi^{\sigma}(h)=
\operatorname{Tr}(\sigma(\xi)\pi(\xi,h))$ associated with
$\sigma$. Then we have by (\ref{eq:wave}) (with $h\in\H$)
\begin{equation}\label{eq:JisW}
\J(\tilde{\sigma})(h)=\mu_{\Ri,\d}W_{\sigma}(h),
\end{equation}
where we put for any
$\sigma\in\cc(\Xi_{P,\d,u},\operatorname{End}(\V_\Xi))$,
\begin{align}\label{eq:wavepack}
\index{W@$W_{\sigma}$, wave packet}
W_{\sigma}(h):&=\int_{\Xi_{P,\d,u}}\Phi_\xi^\sigma(h)c^{-1}(\xi)d\xi\\
\nonumber     &=\int_{\Xi_{P,\d,u}}
\operatorname{Tr}(\sigma(\xi)\pi(\xi,h))c^{-1}(\xi)d\xi
\end{align}
\begin{thm}\label{thm:mainest}
For every $k\in\Z_+$, there exists a continuous seminorm
$p_k$ on $\cc(\Xi_{P,\d,u},\operatorname{End}(\V_\Xi))$
such that
\begin{equation}\label{eq:mainest}
\vert W_{\sigma}(N_{u}\theta_{x}N_{v})\vert \leq
(1+\Vert x\Vert )^{-k}p_{k}(\sigma),
\end{equation}
for all $x\in X^{+}$,
$u,v\in W_{0}$ and
$\sigma\in\cc(\Xi_{P,\d,u},\operatorname{End}(\V_\Xi))$.
\end{thm}
\begin{proof}
First, by using right and left translations by the $N_{w}$, $w\in
W_{0}$,  and Corollary \ref{cor:trans},
it is enough to prove (\ref{eq:mainest}) for $u=v=1$. Thus,
we assume
$u=v=1$ in the following.

The proof is by induction on the rank of $X$. The statement
is clear if the rank of $X$ is zero. One assumes the theorem is
true for lattices of rank strictly smaller than rank of $X$.

For the induction step we consider two cases, namely
the case where $Z_X\not=0$ (first case), and the case
where $Z_X=0$ (second case).

{\it First case.} In this case the semisimple quotient
$\H_{F_0}$ of $\H=\H^{F_0}$ has smaller rank than $\H$.
Recall the results of Proposition \ref{prop:parber} and
Proposition \ref{prop:partwist}.
Let us
denote the semisimple quotient $\H_{F_0}$ of $\H$
by $\H_0$, its root datum $\Ri_{F_0}$ by $\Ri_0$ etc.

We have $T_P\subset T_0$ and $T^P\supset T^0$. Let
$T^P_0=(T_0)^P$ be the connected component of $e$ of the
intersection $T_0\cap T^P$. Then the product map
$T_P\times T_0^P\to T_0$ is a finite covering, as
is the product map $T_0^P\times T^0\to T^P$. Let
$\xi=(P,\d,t^P)$ and suppose that $t^P=t_0^Pt^0$
for $t_0^P\in T^P_{0,u}$ and $t^0\in T^0_u$.
Let $\xi_0=(P,\d,t_0^P)\in \Xi_{\Ri_0,P,\d,u}$
denote the standard induction datum for $\H_0$.
Recall the epimorphism
$\phi_{t^0}:\H\to\H_0$ of Proposition
\ref{prop:partwist}.
It is easy to see that $\pi(\xi)=\pi(\xi_0)_{t^0}$, the
pull back of the representation $\pi(\xi_0)$ of $\H_0$
to $\H$ via $\phi_{t^0}$. This implies that
\begin{equation}
\pi(\xi)(\theta_x)=t^0(x)\pi(\xi_0)(\theta_{x_0})
\end{equation}
for all $x\in X$, where $x_0\in X_0$ is the
canonical image of $x$ in $X_0$.

Hence, since $c(\xi)=c(\xi_0)$
(indeed, use Definition \ref{dfn:cind} and observe that
$\a(t)=\a(r_Pt^P)=\a(r_Pt^P_0t^0)=\a(r_Pt^P_0)$
for all $\a\in R_0$)
and since
\begin{equation}
\int_{T^P_u}f(t^P)dt^P=\int_{T^P_{0,u}\times
T^0_u}f(t^P_0t^0)dt^P_0dt^0
\end{equation}
for all integrable
functions $f$ on $T^P_u$, we have
\begin{equation}\label{eq:Wred}
W_\sigma(\theta_x)=W_{0,\sigma_x}(\theta_{x_0}),
\end{equation}
where
$\sigma_x\in\cc(\Xi_{\Ri_0,P,\d,u},\operatorname{End}(\V_{\Xi_0}))$
is defined by
\begin{align}\label{eq:sx}
\sigma_x(\xi_0)&=\int_{T^0_u}t^0(x)\sigma(t^P_0t^0)dt^0\\
\nonumber &=\int_{T^0_u}t^0(x^0)\sigma(t^P_0t^0)dt^0,
\end{align}
and where $x^0$ is the canonical image of
$x$ in $X^0$.

From equation (\ref{eq:sx}) it is clear, by harmonic analysis on
the torus $T_{0,u}^P\times T^0_u$, that for all $k\in \Z_+$
and all continuous seminorms $q$ on
$\cc(\Xi_{\Ri_0,P,\d,u},\operatorname{End}(\V_{\Xi_0}))$,
there exists a continuous seminorm $p=p_{q,k}$ on
$\cc(\Xi_{P,\d,u},\operatorname{End}(\V_{\Xi}))$
such that
\begin{equation}\label{eq:partfour}
q(\sigma_x)=(1+\Vert x^0\Vert)^{-k}p(\sigma)
\end{equation}
for all $x\in X$ and for all
$\sigma\in\cc(\Xi_{P,\d,u},\operatorname{End}(\V_{\Xi}))$.

Now apply the induction hypothesis to $W_{0,\sigma_x}$.
In view of (\ref{eq:Wred}) and (\ref{eq:partfour}) this
yields the induction
step in the first case.

{\it Second case.} We now consider the case $Z_X=0$.
If $\alpha \in F_{0}$ and $a>0$,
one defines $X^{+}_{\alpha, a}=\{x\in X^{+}\mid
\langle x,\a^\vee\rangle>a\Vert
x\Vert\}$. We first prove that:
\begin{equation}\label{eq:all}
\cup _{\alpha\in F_{0}}X^{+}_{\alpha, a}=
X^{+}\setminus \{0\},\mathrm{\ for\ all\ }a\mathrm{\ small\ enough}
\end{equation}
Let $x\in X^{+}\setminus\{0\}$. We write
$$x= \sum _{\alpha\in F_{0}}\langle x,\a^\vee\rangle\d_{\a}$$
where $(\d_{\a})$ are the fundamental weights. For $a>0$ small
enough, one has, by equivalence of norms in finite dimensional
vector spaces:
\begin{equation}\label{eq:equivnorm}
2a\Vert x\Vert \leq \operatorname{sup}_{\a \in
F_{0}}\vert \langle x,\a^\vee\rangle\vert, \>\> x \in X^{+}
\end{equation}
Then, for
$x\in X^{+}\setminus\{0\}$, choose $\a \in F_{0}$ with
$\langle x,\a^\vee\rangle$
maximal.
From  (\ref{eq:equivnorm}), one has $\langle x,\a^\vee\rangle
\geq 2a \Vert x\Vert$,
hence, as $\Vert x\Vert \not= 0$: $$\langle x,\a^\vee\rangle
> a \Vert x\Vert,
\>\>i.e. \>\>x \in X^{+}_{\alpha, a} $$ which proves
(\ref{eq:all}).

Hence it is enough to prove the estimates for
$x\in  X^{+}_{\a, a}$. Let $Q=F_{0}\setminus\{\a\}$. Then it
follows from Lemma \ref{lem:est1},
that for some $b>0$, and $C>0$, one has:
$$\vert \Phi^\sigma_{t^{P}}( \theta_{x})- \Phi^{\sigma,Q}_{t^{P}}(
\theta_{x})\vert  \leq Ce^{-b\Vert x\Vert }, \>\>
\mathrm{\ for\ all\ }x \in X^{+}_{\a,
a}, \>\>t^{P}\in T^{P}_{u}$$  By integration of this inequality
over $T^P_u$ against the continuous function $|c^{-1}(\xi)|$
it suffices to prove the estimates (\ref{eq:mainest})
after replacing  $\Phi_{t^{P}}^\sigma$ by
$\Phi^{\sigma,Q}_{t^{P}}$. But by Corollary \ref{cor:fina},
the restriction to ${\mathcal H}^{Q}$ of the constant term
$\Phi^{\sigma,Q}_{t^{P}}c^{-1}(\xi)$ of the
normalized smooth
family $\Phi_{t^{P}}^\sigma c^{-1}(\xi)$ of coefficients
is a sum of normalized smooth families of coefficients for
${\mathcal R}_{Q}= (X,Y,R_{Q}, {R}^\vee_{Q},Q)$. This brings us
back at the first case of the induction step, thus finishing the
proof.
\end{proof}
\begin{cor}\label{cor:mainJ}
It follows from Theorem \ref{thm:mainest}, (\ref{eq:JisW})
and Lemma \ref{lem:w=uxv} that $\J(\sigma)\in\S$ for all
$\sigma\in
\mathcal{C}(\Xi_{P,\d,u},\operatorname{End}(\V_\Xi))$,
and that $\J_\mathcal{C}:
\mathcal{C}(\Xi_{P,\d,u},\operatorname{End}(\V_\Xi))\to\S$
is contiunous.

By Lemma \ref{lem:ccc} we see that in particular
$\J(\sigma)\in\S$ for all
$\sigma\in\cc(\Xi_{P,\d,u},\operatorname{End}(\V_\Xi))$.
\end{cor}
\subsection{End of the proof of the main Theorem}
We start with a basic technical lemma:
\begin{lm}\label{lem:w=uxv}
Let $n\in \Z $.  There exists a constant $C_{n}$ with the following
property. For all $f\in {\H}^{*}$ for which there exists $C>0$ such
that:
\begin{equation}\label{eq:assu}
\vert f(T_{u}\theta_{x}T_{v}) \vert \leq C
 (1+\Vert x\Vert )^{-n}, \>\> u,v\in W_{0}, x\in  X^{+}
\end{equation}
then one has:
\begin{equation*}
\vert f(N_{w}) \vert \leq C_{n}C(1+ {\mathcal N}(w))^{-n}, \>\> w\in W
\end{equation*}
\end{lm}
\begin{proof}
As in \cite{O} (2.25), one writes, for $w=uxv$, with $  u,v\in
W_{0}, x\in  X^{+}$,
\begin{equation} N_{w}=\sum_{  u,v\in W_{0}}c_{w, (u', v')}
  T_{u'}\theta_{x}T_{v'}
\end{equation}
where the real  coefficients $c_{uxv, (u', v')}$ and
$c_{uyv, (u', v')}$ are equal if $x$ and $y$ belong to the same facets
of the cone $X^{+}$. The number of facets being finite, one sees, by
using the assumption (\ref{eq:assu}), that there exists $C'$ such that:
\begin{equation}\label{eq:eqs1}
\vert f(N_{uxv}\vert\leq C'C (1+\Vert x\Vert )^{-n},\>\>
  u,v\in W_{0}, x\in  X^{+}
\end{equation}
But, from \cite{O} (2.27), one deduces the existence of
$r_{0}\geq 0$ such that:
\begin{equation}\label{eq:eqs2}
{\mathcal N}(x)-r_{0}\leq {\mathcal N}(uxv)\leq {\mathcal N}(x)+r_{0},
 \>\> u,v\in
W_{0}, x\in  X^{+}
\end{equation}
But one has :
${\mathcal N}(x)=\langle x,2{\rho^\vee}\rangle+
\Vert x^{0}\Vert, \>\> x\in X^{+}$
where $x^{0}$ is the projection of $x\in X\otimes \R$ on
$Z_{X} \otimes \R$ along $\Z R_{0}\otimes \R$. Let us define
$$ \Vert {\mathit v}\Vert^\prime=
\operatorname{sup}_{u \in W_{0}}\vert {\mathit
v}(2{u\rho^\vee})
\vert +\Vert {\mathit v}^{0}\Vert,\>\>  {\mathit v}\in X\otimes \R $$
Then $\Vert\cdot\Vert^\prime$ is a norm on $X\otimes \R$, which is equivalent
to $\Vert\cdot\Vert$. Moreover
$${\mathcal N}(x)=\Vert x\Vert^\prime, \>\> \mathrm{\ for\ all\ }x\in X^{+}.$$
Hence there exists $C''>0$ such that:
\begin{equation}\label{eq:eqs3}
{C''}^{-1}{\mathcal N}(x)\leq \Vert x \Vert \leq C'' {\mathcal N}(x)
\end{equation}
Taking into account (\ref{eq:eqs1}), (\ref{eq:eqs2})  and (\ref{eq:eqs3}),
one gets the result.
\end{proof}
\ { \em End of the proof of the Main Theorem:}\ste
By Corollary \ref{cor:imF}, the image of $\F_\S$ is contained in
the space of smooth $\W$-equivariant sections
$\cc(\Xi_{u},\operatorname{End}(\V_\Xi))^\W$, and $\F_\S$
is continuous.

Corollary \ref{cor:mainJ} states that
the image of $\J_\mathcal{C}$ is contained in $\S$, and
that $\J_\mathcal{C}
:\mathcal{C}(\Xi_{u},\operatorname{End}(\V_\Xi))
\to\S$ is continuous.

Since
$\cc(\Xi_u,\operatorname{End}(\V_\Xi))^\W\subset
\mathcal{C}(\Xi_u,\operatorname{End}(\V_\Xi))$ (see Lemma
\ref{lem:ccc}) and $\S\subset L_2(\H)$ (see (\ref{eq:sinl2})),
Corollary \ref{cor:wave} implies that
$\J_\mathcal{C} \F_\S=\operatorname{id}_\S$.
It follows that
the map $\J_\mathcal{C}$ in (\ref{eq:js}) is surjective, and
that $\F_\S$ is injective.

Since $\mathcal{C}\subset L_2(\Xi_u,\operatorname{End}(\V_\Xi),\mu_{Pl})$
(see Lemma \ref{lem:ccc}),
Corollary \ref{cor:wave} also implies that
$\F_\S \J_\mathcal{C}=p_{\W,\mathcal{C}}$.
It follows, since $p_{\W,\mathcal{C}}$
is the identity on
$\cc(\Xi_u,\operatorname{End}(\V_\Xi))^\W\subset
\mathcal{C}(\Xi_u,\operatorname{End}(\V_\Xi))$, that $\F_\S$
is also surjective in (\ref{eq:fs}).
This finishes the proof of the Main Theorem.
\qed
\section{Appendix: Some applications of spectral projections}
The following lemma was suggested
by Lemma 20.1 from \cite{B} and its proof.
\begin{lem}\label{lem:ban}
Let $V$ be a
complex normed vector space of dimension $p$.  There exists $C>0$ such
that for all $A\in \operatorname{End}(V)$ with eigenvalues
less or equal to $1$:
\begin{equation}
\Vert A^n\Vert\leq C (1+\Vert A\Vert )^{p-1}(1+n)^{p}, n\in
\Z_+
\end{equation}
Here $\Vert A\Vert$ is the operator norm of $A$.
\end{lem}
\begin{proof}
Let $D_{n}$ be the  disk of center $0$ and radius
$1+(1+n)^{-1}$.
Then
\begin{equation}
A^{n}=1/2i\pi\int _{\partial
D_{n}}z^{n}(z\operatorname{Id}-A)^{-1}dz\end{equation}
From the Cramer rules, there
exists a polynomial function from $\operatorname{End}(V)$
into itself, $B\mapsto M(B)$,
of degree $p-1$, such that for any invertible $B$, one has:
\begin{equation}B^{-1}=(\operatorname{det}(B))^{-1}M(B)\end{equation}
Hence, there
exists $C'>0$ such that: $$\Vert M(B)\Vert\leq C' (1+\Vert B\Vert
)^{p-1}, B\in \operatorname{End}(V).$$
Hence, taking into account: $$1+\Vert z\operatorname{Id}-A\Vert
\leq 2+(1+n)^{-1}+\Vert A\Vert\leq (2+(1+n)^{-1})(1+\Vert A\Vert), z\in
D_{n},  $$ one has
\begin{equation}  \Vert M(z\operatorname{Id}-A)\Vert \leq C'
(2+(1+n)^{-1})^{p-1}(1+\Vert A\Vert)^{p-1}, z\in D_{n}
\end{equation}
Now the eigenvalues of $z\operatorname{Id} -A$,
$z\in \partial D_{n}$ are of modulus
greater or equal to  $(1+n)^{-1}$. Hence
\begin{equation}\vert
\operatorname{det}(z\operatorname{Id}-A)\vert \geq (1+n)^{-p}, n\in
\Z_+
\end{equation}
The length of
$\partial D_{n}$ is $2\pi(1+(1+n)^{-1})$.  From equations (6.1) to
(6.4),  one gets:
$$\Vert A^n\Vert \leq (1+(1+n)^{-1})^{1+n}
C'(2+(1+n)^{-1})^{p-1}(1+n)^{p}(1+\Vert A\Vert)^{p-1}$$
One gets the required
estimate with: $$C=C' e3^{p-1}$$
\end{proof}
\begin{cor}\label{cor:bancor}
\begin{enumerate}
\item[(i)] Let
$r'>r>0$. There exists $C_{r,r'}$ such that for all
$A\in \operatorname{End}(V)$
with eigenvalues of modulus less or equal to $r$, one has: $$\Vert
A^{n}\Vert\leq C_{r,r'}{(r')}^{n}(1+\Vert A\Vert )^{p-1}$$
\item[(ii)]
Let $\epsilon >0$ and let $\Omega_{\epsilon }$ be the set of
elements in $\operatorname{End}(V)$
whose eigenvalues are either of modulus one or
of modulus less or equal to $1-\e$. Let $P_{<1}$ be the sum of the
spectral projections corresponding to the eigenvalues of modulus
strictly less than 1. Then $P_{<1}A^{n}= (A_{<1})^{n}$, where
$A_{<1}=P_{<1}A$. Let $b>0$ such that $1-\epsilon< e^{-b}$. There
exists $C$ depending on $ b, \e$ and $V$ such that:  $$ \Vert
P_{<1}A^{n}\Vert\leq C( 1+\Vert P_{<1}A\Vert) ^{p-1}e^{-bn}, n\in
\N, A\in \Omega_{\e}$$
\item[(iii)]
If A(t) is a
continuous ( resp. holomorphic) function with values  in
$\Omega_{\e}$, then $A_{<1}(t)$ has the same  property
\end{enumerate}
\end{cor}
\begin{proof}
(i) One applies the previous result to ${r''}^{-1}A$, where
$r<r''<r'$
and one uses the fact that $(1+n)^{p}(r'/r'')^{-n}$ is bounded.

(ii) follows from (i) applied to $A_{<1}$, $r=1-\epsilon$, $r'=
e^{-b}$.

(iii) follows from the formula
$$ A_{<1}= 1/2i\pi \int _{\partial D}z(z\operatorname{Id}-A)^{-1}dz$$
where $D$ is the disc of center 0 and radius $1-\e/2$.
\end{proof}
\begin{lem}\label{lem:specproj}
Let $\epsilon, a>0$, $p, l\in \N$. Let $V$ a normed complex vector
space of dimension $p$ and $X=\Z^{l}$. Let $\pi$ be a finite
dimensional complex  representation of $X$. One denotes by
$(e_{1},\dots,
e_{p})$ the canonical basis of $X$.  One sets  $A_{1}=\pi(e_{1}),
\dots, A_{l}=\pi(e_{l})$. If $n=(n_{1},\dots , n_{l})\in X$, one
sets :
$\Vert n\Vert = \vert n_{1}\vert+\dots +\vert n_{l}\vert$, and
$A^{n}=\pi(n)$.\ste
Assume that the modulus of the eigenvalues of
the $A_{i}$ are  less or equal to one, and the eigenvalues of
$A_{1}$ are either of modulus one or of modulus less or equal to
$1-\epsilon$. Let us denote by $P_{>1}$ (resp. $P_{1}$) the sum
of the spectral projections of $A_{1}$ corresponding to the
eigenvalues of modulus strictly less than 1 (resp. of modulus
1).

Set $X^{+}_{a}=\{ n \in \Z_+^{l}\vert  n_{1}> a \Vert
n\Vert\}$.\ste  Then there exists $a'$ and $C'$,
independent of the
representation $\pi$ in $V$, such that:
$$\Vert P_{<1}A^{n}\Vert \leq C'(\prod_{i=1, \dots, l}(1+\Vert
A_{i}\Vert )^{p})e^{-a' \Vert n\Vert}, \>\> n\in X^{+}_{a}$$
\end{lem}
\begin{proof}
 From Corollary \ref{cor:bancor}(i), (ii), one deduces that, for
$b>0$ such that $1-\epsilon< e^{-b}$, and $b'>0$, there exists a
constant $C>0$, depending only on $\e, b, b'$ and $V$ such that:
$$\Vert P_{<1}A^{n}\Vert \leq C'(\prod_{i=1, \dots, l}(1+\Vert
A_{i}\Vert )^{p})e^{-bn_{1}+b'(n_{2}+\dots +n_{l})}, n \in X^{+}$$
If $n\in X^{+}_{a}$, one has :
$$bn_{1}-b'(n_{2}+\dots +n_{l})\geq (ab-b')\Vert n\Vert$$
$b$ being chosen, one takes $b'= ab/2$. Then the inequality of the
Lemma is satisfied for $a'=ab/2$.
\end{proof}
\section{Appendix: The $c$-function}
In this appendix we have collected some of the properties of the
Macdonald $c$-function. These properties play a predominant
role throughout this paper, and are closely related with the
properties of residual cosets as discussed in \cite{O},
Appendix: residual cosets (Section 7).

The Macdonald $c$-function is defined as
the following expression
$c\in{}_\mathcal{Q}\A=
\mathcal{Q}\otimes_\Z\A$\index{A@${}_\mathcal{Q}\A=
\mathcal{Q}\otimes_\Z\A$, quotient field of $\A$}
the quotient field of $\A$:
\begin{equation}\label{eq:defc}
\index{c@$c(\xi)=
\prod_{\a\in R_{0,+}\backslash R_{P,+}}c_\a(t)$,
Macdonald's $c$-function on $\Xi$, where $\xi=(P,\d,t^P)$,
and $t=r_Pt^P$ with $W_Pr_P$ the
central character of $\d$}
c:=\prod_{\a\in R_{0,+}}c_\a=\prod_{\alpha\in R_{1,+}}c_\a,
\end{equation}
where $c_\a$\index{c@$c_\a$, rank one $c$-function}
for $\a\in R_1$ is equal to
\begin{equation}\label{eq:defcr1}
\index{c@$c_\a$, rank one $c$-function}
c_\a:= \frac{(1+q_{\alpha^\vee}^{-1/2}\theta_{-\alpha/2})
(1-q_{\alpha^\vee}^{-1/2}q_{2\alpha^\vee}^{-1}\theta_{-\alpha/2})}
{1-\theta_{-\alpha}}\in{}_\mathcal{Q}\A.
\end{equation}
If $\a\in R_0\backslash R_1$ then we define
$c_\a:=c_{2\a}$.
\begin{rem}\label{rem:conv}
We have thus associated a $c$-function $c_\a$ to each root
$\a\in R_{nr}$, but $c_\a$ only depends on the direction of
$\a$. This convention was used in \cite{O}, but differs from the
one used in \cite{EO}. If $\a\in R_1$ and $\a/2\not\in R_0$, then
the formula for $c_\a$ should be interpreted by setting
$q_{2\alpha^\vee}=1$,
and then rewriting the numerator as
$(1-q_{\a^\vee}^{-1}\theta_{-\alpha})$.
Here and below we use this convention.
\end{rem}
We view $c$ as a rational function on $T$ via the isomorphism
of $\A$ and the ring of regular functions on $T$ sending
$\theta_x$ to the complex character $x$ of $T$.

Since the numerator and the denominator of $c$ both are
products of irreducible factors whose zero locus is nonsingular
(a coset of a codimension $1$ subtorus of $T$), it is straightforward
to define the pole order
$i_t$\index{i@$i_t$, pole order of $(c(t)c(w_0t))^{-1}$
at $t\in T$}
of $(c(t)c(w_0t))^{-1}$
at a point $t\in T$ (see \cite{O}, Definition 3.2).

Let $Q=Q(R_0)$ denote the root lattice of $R_0$. The following theorem
is the main property of the $c$-function:
\begin{thm}\label{thm:mainc}(\cite{O}, Theorem 7.10)
We have: $i_t\leq\operatorname{rank}(Q)$ for all $t\in T$.
\end{thm}
We define the notion of a ``residual point'' of $T$ (with
respect to $(\Ri,q)$) as follows:
\begin{dfn}\label{dfn:res}
A point $t\in T$ is called residual if $i_t=\operatorname{rank}(X)$.
\end{dfn}
\begin{cor}\label{cor:finc}
There exist only finitely many residual points in $T$, for every
root datum $\Ri$ and label function $q$ for $\Ri$, and the set of
residual points is empty unless $Z_X=0$.
\end{cor}
\begin{proof}
Let $n=\operatorname{rank}(R_0)$.
From equation (\ref{eq:defc}) it is clear that for any $k\in\Z$,
the set $S_k:=\{t\mid i_t=k\}$ is a finite (possibly empty)
union of nonempty Zariski
open subsets of cosets $L$ of subtori of $T$, whose Lie algebra is an
intersection of root hyperplanes $\a=0$ of $\mathbb{C}\otimes Y$.
If $L$ is such a coset with $\operatorname{codim}(L)=d$, then
$R_L:=\{\a\in R_0\mid \a|_L{\rm\ is\ constant}\}$ is a parabolic
subsystem of rank $d$. Moreover, the projection $t_L$ of $L$
onto $T_L$ is point with $i^{R_L}_{t_L}=k$. Applying Theorem
\ref{thm:mainc} to $T_L$ with respect to $(\Ri_L,q_L)$ we obtain
$k\leq d$.

Hence if $S_n$ is not finite, then there exists a proper parabolic
root subsystem $R_L\subset R_0$ of rank $m<n$ say, such that $n\leq
m$, which is clearly absurd.
The remaining part of the Corollary is straightforward from
Theorem \ref{thm:mainc}.
\end{proof}
\begin{rem}
There is a classification of residual points (cf. \cite{HOH0}, \cite{O}).
\end{rem}
Another fact of great consequence is the following.
\begin{thm}\label{thm:cster}(\cite{O}, Theorem 7.14)
Let $r=sc\in T$ be residual, with $s\in T_u$ and $c\in T_{rs}$.
Then $r^*:=sc^{-1}\in W(R_{s,1})r$, where
$R_{s,1}$ is the root subsystem of $R_1$ defined by
$R_{s,1}:=\{\a\in R_1\mid \a(s)=1\}$.
\end{thm}
We extend the definition of the $c$-function to arbitrary
standard induction data. First recall Theorem \ref{thm:dsres},
stating that the central character of a residual discrete series
representation is residual.
\begin{dfn}\label{dfn:cind}
Let $\xi=(P,\d,t^P)$ be a standard induction datum, and
let $r_P\in T_P$ be such that $W_Pr_P$ is the central
character of $\d$ (thus $r_P$ is a $(\Ri_P,T_P)$-residual point).
Put $t=r_Pt^P\in T$. We define:
\begin{equation}\label{dfn:c2}
c(\xi):=\prod_{\a\in R_{0,+}\backslash R_{P,+}}c_\a(t)
\end{equation}
Notice that we recover the original $c$-function defined
on $T$ as the special case where $P=\emptyset$ and $\d=\C$.
\end{dfn}

The next proposition goes back to \cite{HOH0}, Theorem 3.13
(also see \cite{O}, Theorem 3.25).
\begin{prop}\label{lem:csmo}
Let $P\subset F_0$ and
let $\xi=(P,\d,t^P)\in\Xi_{P,\d,u}$. Choose $r_P\in T_P$
such that $W_Pr_P$ is the central character of $\d$, and
let $t=r_Pt^P\in T$.
\begin{enumerate}
\item $c(\xi^{-1})=c(w^P(\xi))=\overline{c(\xi)}$.
\item The function $\xi\to \vert c(\xi)\vert^{2}$ on $\Xi_u$
is $\W$-invariant.
\item The function $c(\xi)$ is $\K$-invariant.
\item Let $P^\prime\in\P$ and
$d\in K_{P^\prime}\times W(P,P^\prime)$.
The rational functions $c(d\xi){c(\xi)}^{-1}$
and $c(d\xi)^{-1}{c(\xi)}$ (of $\xi\in\Xi_{P,\d}$) are regular in
a neighborhood of $\Xi_{P,\d,u}$.
\item The rational function ${c(\xi)}^{-1}$ is regular in
a neighborhood of $\Xi_{u}$.
\end{enumerate}
\end{prop}
\begin{proof} (i) A straightforward computation from
the definitions, using Theorem \ref{thm:cster}
(cf. \cite{O}, (3.58)).

(ii) The $\W$-invariance follows simply
from the definitions if we write (using (i))
$\vert c(\xi)\vert^{-2}=(c(\xi)c(\xi^{-1})^{-1}$
(cf. \cite{O}, Proposition 3.27).

(iii) This follows trivially from the definition of the
action of $\K$ on $\xi$: If $k\in K_P$ then
$k\xi=k(P,\d,t^P)=(P,\Psi_k(\d),kt^P)$. The central character
of $\Psi_k(\d)$ is equal to $k^{-1}W_Pr_P=W_P(k^{-1}r_P)$, thus
we need to evaluate the $c_\a$ in the product
$c(k\xi)$ at the point $k^{-1}r_Pkt^P=t$, or any of its
images under the action of $W_P$. Hence $c(k\xi)=c(\xi)$.

(iv) By (i) and (ii) it is clear that these
rational functions have modulus $1$ on $\Xi_{P,\d,u}$
(outside their respective singular sets).

The singular sets of these rational functions are of the
following form.
Choose $r_P\in T_P$ such that $W_Pr_P$ is the central
character of $\d$.
Then the singular set of $c(d\xi){c(\xi)}^{-1}$ is the union
of the zero sets of the functions
\begin{equation}
\prod_{\a\in R_{1,+}\backslash R_{P,1,+}}(1-\a_{P,\d}^{-1})
\end{equation}
and
\begin{equation}
\prod_{\a\in R_{1,+}\backslash R_{P,1,+}}
(1+q_{\a^\vee}^{-1/2}\a_{P,\d}^{-1/2})
(1-q_{\a^\vee}^{-1/2}q_{2\a^\vee}^{-1}\a_{P,\d}^{-1/2})
\end{equation}
on $\Xi_{P,\d}$,
where $\a_{P,\d}$ denotes the function on $\Xi_{P,\d}$
defined by $\a_{P,\d}(\xi)=\a(r_Pt^P)$.

In the case of $c(w^P\xi)^{-1}{c(\xi)}$ the answer is the
same, but we need to take the products over the set
$\a\in d^{-1}R_{1,+}\backslash R_{P,1,-}$.

The intersection of a component of this hypersurface
with $\Xi_{P,\d,u}$ is either empty or it has (real) codimension $1$
in $\Xi_{P,\d,u}$.

By the boundedness of $c(d\xi){c(\xi)}^{-1}$ on $\Xi_{P,\d,u}$,
this implies that the pole order of this function at a
components of the singular set which meets $\Xi_{P,\d,u}$ is in fact
equal to zero. Hence the poles are removable in a neighborhood of
$\Xi_{P,\d,u}$. Similarly for $c(d\xi)^{-1}{c(\xi)}$.

(v) The proof of \cite{HOH0}, Theorem 3.13 may be adapted to
the present situation. Or we may argue as in (iv) as follows.

Since $|c(\xi)|^{-2}=(c(\xi)c(w^P(\xi)))^{-1}$ is smooth
on $\Xi_u$ (cf. \cite{O}, Theorem 3.25, equation (3.53),
Proposition 3.27 and equation (3.58)),
it follows that $c(\xi)^{-1}$ is bounded on $\Xi_{P,\d,u}$.
Hence the argument that was used in the proof of (iv) applies.
\end{proof}

\printindex
\end{document}